\documentclass[12pt]{amsart}

\usepackage{amscd,amssymb}
\usepackage{amsthm,amsmath,amssymb}

\newtheorem{theorem}[equation]{Theorem}

\newtheorem{proposition}[equation]{Proposition}
\newtheorem{lemma}[equation]{Lemma}
\newtheorem{corollary}[equation]{Corollary}
\newtheorem{conjecture}[equation]{Conjecture}

\theoremstyle{definition}
\newtheorem{example}[equation]{Example}
\newtheorem{definition}[equation]{Definition}

\theoremstyle{remark}
\newtheorem{remark}[equation]{Remark}

\author{Ivan Cheltsov}

\title{On factoriality of nodal threefolds}

\address{Steklov Institute of Mathematics\hfill\break\indent 8 Gubkin street, Moscow 117966\hfill\break\indent Russia\hfill}%
\address{The University of Edinburgh\hfill\break\indent School of Mathematics\hfill\break\indent Kings Buildings, Mayfield Road\hfill\break\indent Edinburgh EH9 3JZ, UK\hfill\break\indent}%

\email{cheltsov@yahoo.com, I.Cheltsov@ed.ac.uk}

\thanks{The author is very grateful to A.Corti, S.\,Cynk, M.Grinenko,
V.Is\-kov\-skikh, S.Kudryavtsev, V.Kulikov, M.Mel\-la, J.Park,
Yu.Pro\-kho\-rov, A.Pukh\-li\-kov, V.Sho\-ku\-rov and L.Wotzlaw
for fruitful conversations. A special thank goes to V.Kulikov for
Lemmas~\ref{lemma:double-solid-III} and
\ref{lemma:hypersurface-III}.}

\thanks{All varieties are assumed to be
projective, normal, and defined over $\mathbb{C}$.}

\begin{document}

\begin{abstract}
We prove the $\mathbb{Q}$-factoriality of a nodal hypersurface in
$\mathbb{P}^{4}$ of degree $n$ with at most
${\frac{(n-1)^{2}}{4}}$ nodes and  the $\mathbb{Q}$-factoriality
of a double cover of $\mathbb{P}^{3}$ branched over a nodal
surface of degree $2r$ with at most ${\frac{(2r-1)r}{3}}$ nodes.
\end{abstract}

\maketitle


\section{Introduction.}
\label{section:introduction}

Nodal 3-folds\footnote{A 3-fold is called nodal if all its
singular points are ordinary double points.} arise naturally in
many different topics of algebraic geometry. For example, the
non-rationality of many smooth rationally connected 3-folds  (see
\cite{Ko96}) is proved via the degeneration to nodal 3-folds (see
\cite{Be77}, \cite{Tu80}, \cite{Ch04c}, \cite{Ch04a}). Obviously,
nodal 3-folds are the simplest degenerations of smooth ones.
However, the geometry can be very different in smooth and nodal
cases.

Every surface in a smooth hypersurface in $\mathbb{P}^{4}$ is a
complete intersection by the Lefschetz theorem (see
\cite{AndFr59}, \cite{Bott59}, \cite{Vo03}), which is not the case
if the hypersurface is nodal. The birational automorphisms of a
smooth quartic 3-fold in $\mathbb{P}^{4}$ form a finite group
consisting of projective automorphisms (see \cite{IsMa71}), but
for any non-smooth nodal quartic 3-fold this group is always
infinite (see \cite{Pu88b}, \cite{Co00}, \cite{Me03}). Every
smooth quartic 3-fold and every smooth cubic 3-fold in
$\mathbb{P}^{4}$ are non-rational (see \cite{IsMa71},
\cite{ClGr72}), but all singular nodal cubic 3-folds are rational
and there are rational nodal quartic 3-folds (see \cite{Pet98}).
Every smooth double cover of $\mathbb{P}^{3}$ ramified in a sextic
or quartic surface is non-rational (see \cite{Is80}, \cite{Ti80b},
\cite{Ti82}, \cite{Ti86}, \cite{Cl91}), but there are rational
nodal ones (see \cite{Kr00}, \cite{ChPa04}).

The simplest examples of nodal 3-folds are nodal hypersurfaces in
$\mathbb{P}^{4}$ and double covers of $\mathbb{P}^{3}$ branched
over a nodal surfaces. The latter are called double solids. These
3-folds were studied in \cite{Cl83}, \cite{Sch85}, \cite{We87},
\cite{Fi87}, \cite{FiWe89}, \cite{Di90}, \cite{Bo90},
\cite{vSt93}, \cite{Pet98}, \cite{Cy99}, \cite{Cy01}, \cite{Cy02},
\cite{Me03}, \cite{CiGe03}, \cite{ChPa04}.

For a given nodal 3-fold, it is one of substantial questions
whether it is $\mathbb{Q}$-factorial\footnote{A variety is called
$\mathbb{Q}$-factorial if a multiple of every Weil divisor on the
variety is a Cartier divisor.} or not. This global topological
property has very simple geometrical description. Namely, a
three-dimensional ordinary double point admits two small
resolutions that differs by a simple flop (see \cite{Ka88},
\cite{We87}, \cite{Ko89}). In particular, a nodal 3-fold with $k$
nodes has $2^{k}$ small resolutions. Therefore the
$\mathbb{Q}$-factoriality of a nodal 3-fold implies that it has no
projective small resolutions.

\begin{remark}
\label{remark:rationality}%
The $\mathbb{Q}$-factoriality of a nodal 3-fold imposes a very
strong geometrical restriction on its birational geometry. For
example, $\mathbb{Q}$-factorial nodal quartic 3-folds and nodal
sextic double solids are non-rational (see \cite{Me03},
\cite{ChPa04}). On the other hand, there are rational
non-$\mathbb{Q}$-factorial nodal quartic 3-folds and nodal sextic
double solids (see \cite{Pet98}, \cite{En99}, \cite{ChPa04}).
\end{remark}

Consider a double cover $\pi:X\to\mathbb{P}^{3}$ branched over a
nodal hypersurface $S\subset\mathbb{P}^{3}$ of degree $2r$ and a
nodal hypersurface $V\subset\mathbb{P}^{4}$ of degree $n$. The
proof of the following result is due to \cite{Cl83}, \cite{We87},
\cite{Di90}, \cite{Cy01}, \cite{Cy02}.

\begin{proposition}
\label{proposition:defect}%
The 3-folds $X$ and $V$ are  $\mathbb{Q}$-factorial if and only if
the nodes of $S\subset\mathbb{P}^{3}$ and $V\subset\mathbb{P}^{4}$
impose independent linear conditions on homogeneous forms of
degree $3r-4$ and $2n-5$ respectively.
\end{proposition}

In particular, $X$ and $V$ are  $\mathbb{Q}$-factorial if
$|\mathrm{Sing}(X)|\le 3r-3$ and $|\mathrm{Sing}(V)|\le 2n-4$
respectively. The main purpose of this paper is to prove the
following two results.

\begin{theorem}
\label{theorem:first}%
Suppose that $|\mathrm{Sing}(X)|\le{\frac{(2r-1)r}{3}}$. Then
$\mathrm{Cl}(X)$ and $\mathrm{Pic}(X)$ are generated by the class
of $\pi^{*}(H)$, where $H$ is a hyperplane in $\mathbb{P}^{3}$. In
particular, $X$ is $\mathbb{Q}$-factorial.
\end{theorem}

\begin{theorem}
\label{theorem:second}%
Suppose that $|\mathrm{Sing}(V)|\le{\frac{(n-1)^{2}}{4}}$. Then
$\mathrm{Cl}(V)$ and $\mathrm{Pic}(V)$ are generated by the class
of a hyperplane section. In particular, $V$ is
$\mathbb{Q}$-factorial.
\end{theorem}

\begin{remark}
\label{remark:factoriality}%
The statements of Theorems~\ref{theorem:first} and
\ref{theorem:second} are equivalent to the
$\mathbb{Q}$-factoriality of the 3-folds $X$ and $V$ respectively.
Indeed, the $\mathbb{Q}$-factoriality of $X$ and $V$ implies
$$
\mathrm{Cl}(X)\otimes\mathbb{Q}\cong\mathrm{Pic}(X)\otimes\mathbb{Q}\cong\mathrm{Cl}(V)\otimes\mathbb{Q}\cong\mathrm{Pic}(V)\otimes\mathbb{Q}\cong\mathbb{Q}
$$
due to the Lefschetz theorem and \cite{Cl83}. Moreover, the groups
$\mathrm{Pic}(X)$ and $\mathrm{Pic}(V)$ have no torsion due to the
Lefschetz theorem and \cite{Cl83}. On the other hand, the local
class group of an ordinary double point is $\mathbb{Z}$ (see
\cite{Mi68}). Therefore, the groups $\mathrm{Cl}(X)$ and
$\mathrm{Cl}(V)$ have no torsion as well (cf. \cite{Cl83}), which
implies the equivalences.
\end{remark}

Actually, the bounds for nodes in Theorems~\ref{theorem:first} and
\ref{theorem:second} are not sharp. For example, in the case $r=3$
the 3-fold $X$ is $\mathbb{Q}$-factorial if $|\mathrm{Sing}(X)|\le
14$ due to \cite{ChPa04}, and in the case $n=4$ the 3-fold $V$ is
$\mathbb{Q}$-factorial if $|\mathrm{Sing}(V)|\le 8$ due to
\cite{Ch04a}.

\begin{example}
\label{example:non-Q-factoriality-I} Consider a hypersurface
$X\subset\mathbb{P}(1^{4},r)$ given by the equation
$$
u^{2}=g^{2}_{r}(x,y,z,t)+h_{1}(x,y,z,t)f_{2r-1}(x,y,z,t)\subset
\mathbb{P}(1^{4},r)\cong\mathrm{Proj}(\mathbb{C}[x,y,z,t,u]),
$$
where $g_{i}$, $h_{i}$, and $f_{i}$ are sufficiently general
polynomials of degree $i$. Let $\pi:X\to\mathbb{P}^3$ be a
restriction of the natural projection
$\mathbb{P}(1^{4},r)\dashrightarrow\mathbb{P}^{3}$, induced by an
embedding of the graded algebras
$\mathbb{C}[x_{0},\ldots,x_{2n}]\subset\mathbb{C}[x_{0},\ldots,x_{2n},y]$.
Then $\pi:X\to\mathbb{P}^3$ is a double cover branched over a
nodal hypersurface
$$g^{2}_{r}(x,y,z,t)+h_{1}(x,y,z,t)f_{2r-1}(x,y,z,t)=0$$
of degree $2r$ and $|\mathrm{Sing}(X)|=(2r-1)r$. Moreover, the
3-fold $X$ is not $\mathbb{Q}$-factorial, because the divisor
$h_{1}=0$ splits into two non-$\mathbb{Q}$-Cartier divisors.
\end{example}

\begin{example}
\label{example:non-Q-factoriality-II} Let $V\subset\mathbb{P}^{4}$
be a hypersurface
$$
xg_{n-1}(x,y,z,t,w)+yf_{n-1}(x,y,z,t,w)\subset
\mathbb{P}^{4}\cong\mathrm{Proj}(\mathbb{C}[x,y,z,t,w]),
$$
where $g_{n-1}$ and $f_{n-1}$ are general polynomials of degree
$n-1$. Then $V$ is nodal and contains the plane $x=y=0$. Hence,
the 3-fold $V$ is not $\mathbb{Q}$-factorial and
$|\mathrm{Sing}(V)|=(n-1)^{2}$.
\end{example}

Therefore, asymptotically the bounds for nodes in
Theorems~\ref{theorem:first} and \ref{theorem:second} are not very
far from being sharp. On the other hand, the following result is
proved in \cite{CiGe03}.

\begin{proposition}
\label{proposition:factoriality} Every smooth surface on $V$ is a
Cartier divisor if $\mathrm{Sing}(V)<(n-1)^2$.
\end{proposition}

Hence, one can expect the following to be true.

\begin{conjecture}
\label{conjecture:factoriality} The inequalities
$|\mathrm{Sing}(X)|<(2r-1)r$ and $|\mathrm{Sing}(V)|<(n-1)^{2}$
imply the $\mathbb{Q}$-factoriality of the 3-folds $X$ and $V$
respectively.
\end{conjecture}

The claim of Conjecture~\ref{conjecture:factoriality} is proved
for $r\le 3$ and $n\le 4$ (see \cite{FiWe89}, \cite{ChPa04},
\cite{Ch04a}). Unfortunately, we are unable to prove
Conjecture~\ref{conjecture:factoriality} in any other case.
However, for every given $r$ and $n$ we always can slightly
improve the bounds in Theorem~\ref{theorem:first} and
\ref{theorem:second}. For example, we prove the following result.

\begin{proposition}
\label{proposition:Calabi-Yau} Let $r=4$ and $n=5$, i.e. $X$ and
$V$ are nodal Calabi-Yau 3-folds, and suppose that
$|\mathrm{Sing}(X)|\le 25$ and $|\mathrm{Sing}(V)|\le 14$. Then
$X$ and $V$ are $\mathbb{Q}$-factorial.
\end{proposition}

The following result is proved in \cite{CiGe03}.

\begin{theorem}
\label{theorem:cili}%
Suppose that the subset $\mathrm{Sing}(V)\subset\mathbb{P}^{4}$ is
a set-theoretic intersection of hypersurfaces of degree
$l<{\frac{n}{2}}$ and $|\mathrm{Sing}(V)|<{\frac
{(n-2l)(n-1)^{2}}{n}}$. Then $V$ is ${\mathbb Q}$-factorial.
\end{theorem}

The saturated ideal of a set of $k$ points in general position in
$\mathbb{P}^{4}$ is generated by polynomials of degree at most
${\frac{n}{4}}$ when $k<(n-1)^{2}$ and $n>72$ by \cite{GeMa84}.
Therefore, Theorem~\ref{theorem:cili} implies the
$\mathbb{Q}$-factoriality of the 3-fold $V$ having less than
${\frac{1}{2}}(n-1)^{2}$ nodes in additional assumption that the
nodes are in general position in $\mathbb{P}^{4}$. However, the
latter generality condition implicitly assumes that the nodes of
$V$ impose independent linear conditions on homogeneous forms of
degree $2n-5$ (see \cite{GeMa84}), which implies the
$\mathbb{Q}$-factoriality of the 3-fold $V$ due to
Proposition~\ref{proposition:defect}. We prove the following
generalization of Theorem~\ref{theorem:cili}.

\begin{theorem}
\label{theorem:six}%
Let $\mathcal{H}\subset|\mathcal{O}_{\mathbb{P}^{3}}(k)|$ and
$\mathcal{D}\subset|\mathcal{O}_{\mathbb{P}^{4}}(l)|$ be linear
subsystems of hypersurfaces vanishing at $\mathrm{Sing}(S)$ and
$\mathrm{Sing}(V)$ respectively. Put
$\hat{\mathcal{H}}=\mathcal{H}\vert_{S}$ and
$\hat{\mathcal{D}}=\mathcal{D}\vert_{V}$. Suppose that
inequalities $k<r$ and $l<{\frac{n}{2}}$ hold. Then
$\mathrm{dim}(\mathrm{Bs}(\hat{\mathcal{H}}))=0$ implies the
$\mathbb{Q}$-factoriality of the 3-fold $X$, and
$\mathrm{dim}(\mathrm{Bs}(\hat{\mathcal{D}}))=0$ implies the
$\mathbb{Q}$-factoriality of the 3-fold $V$.
\end{theorem}

\begin{corollary}
\label{corollary:theorem-six}%
Suppose $\mathrm{Sing}(S)\subset\mathbb{P}^{3}$ and
$\mathrm{Sing}(V)\subset\mathbb{P}^{4}$ are set-theoretic
intersections of hypersurfaces of degree $k<r$ and
$l<{\frac{n}{2}}$ respectively. Then $X$ and $V$ are ${\mathbb
Q}$-factorial.
\end{corollary}

From the point of view of birational geometry the most important
application of Theorems~\ref{theorem:first} and
\ref{theorem:second} is the $\mathbb{Q}$-factoriality condition
for a nodal quartic 3-fold and a sextic double solid, i.e. the
cases $r=3$ and $n=4$ respectively, because in these cases the
$\mathbb{Q}$-factoriality implies the non-rationality (see
\cite{Me03}, \cite{ChPa04}). However, it is possible to apply
Theorems~\ref{theorem:first} and \ref{theorem:second} to certain
higher-dimensional problems in birational algebraic geometry.

\begin{theorem}
\label{theorem:third}%
Let $\tau:U\to\mathbb{P}^{s}$ be a double cover branched over a
hypersurface $F$ of degree $2r$ and $D$ be a hyperplane in
$\mathbb{P}^{s}$ such that $D_{1}\cap\cdots\cap D_{s-3}$ is a
$\mathbb{Q}$-factorial nodal 3-fold, where $D_{i}$ is a general
divisor in $|\tau^{*}(D)|$. Then $\mathrm{Cl}(U)$ and
$\mathrm{Pic}(U)$ are generated by $\tau^{*}(D)$.
\end{theorem}

\begin{theorem}
\label{theorem:forth}%
Let $W\subset\mathbb{P}^{r}$ be a hypersurface of degree $n$ such
that $H_{1}\cap\cdots\cap H_{r-4}$ is a $\mathbb{Q}$-factorial
nodal 3-fold, where $H_{i}$ is a general enough hyperplane section
of $W$. Then the groups $\mathrm{Cl}(W)$ and $\mathrm{Pic}(W)$ are
generated by the class of a hyperplane section of
$W\subset\mathbb{P}^{r}$.
\end{theorem}

A priori Theorems~\ref{theorem:third} and \ref{theorem:forth} can
be used to prove the non-rationality of certain singular
hypersurfaces of degree $r$ in $\mathbb{P}^{r}$ and double covers
of $\mathbb{P}^{s}$ branched over singular hypersurfaces of degree
$2s$ (see \cite{Pu87}, \cite{Pu88a}, \cite{Pu97}, \cite{Ch00b},
\cite{Pu02}, \cite{dFEM03}, \cite{Ch04b}). However, in the former
case the problem can be very hard in general, but in the latter
case the application of Theorems~\ref{theorem:third} can be very
effective. For example, we prove the following result.

\begin{proposition}
\label{proposition:fourfold-double-cover}%
Let $\xi:Y\to \mathbb{P}^{4}$ be a double cover branched over a
hypersurface $F\subset\mathbb{P}^{4}$ of degree $8$ such that $F$
is smooth outside of a smooth curve $C\subset F$, the singularity
of the hypersurface $F$ in sufficiently general point of $C$ is
locally isomorphic to the singularity
$$x_{1}^{2}+x_{2}^{2}+x_{3}^{2}=0\subset\mathbb{C}^{4}\cong\mathrm{Spec}(\mathbb{C}[x_{1},x_{2},x_{3},x_{4}]),$$
the singularities of $F$ in other points of $C$ are locally
isomorphic to the singularity
$$x_{1}^{2}+x_{2}^{2}+x_{3}^{2}x_{4}=0\subset\mathbb{C}^{4}\cong\mathrm{Spec}(\mathbb{C}[x_{1},x_{2},x_{3},x_{4}]),$$
and a general 3-fold in the linear system $|-K_{Y}|$ is
$\mathbb{Q}$-factorial. Then $Y$ is a birationally
rigid\footnote{Namely, the 4-fold $Y$ is a unique Mori fibration
birational to $Y$ (see \cite{Co00}).} terminal
$\mathbb{Q}$-factorial Fano 4-fold with
$\mathrm{Pic}(Y)\cong\mathbb{Z}$ and $\mathrm{Bir}(Y)$ is a finite
group consisting of biregular automorphisms. In particular, the
4-fold $Y$ is non-rational.
\end{proposition}

\begin{example}
\label{example:fourfold} Let $Y\subset\mathbb{P}(1^{5},4)$ be a
hypersurface
$$
u^{2}=\sum_{i=1}^{3}f_{i}(x,y,z,t,w)g^{2}_{i}(x,y,z,t,w)\subset\mathbb{P}(1^{5},4)\cong\mathrm{Proj}(\mathbb{C}[x,y,z,t,w,u]),
$$
where $f_{i}$ and $g_{i}$ are sufficiently general non-constant
homogeneous polynomials such that
$\mathrm{deg}(f_{i})+2\mathrm{deg}(g_{i})=8$. Then the natural
projection $\mathbb{P}(1^{5},4)\dashrightarrow\mathbb{P}^{4}$
induces a double cover $\tau:Y\to\mathbb{P}^4$ branched over a
hypersurface $F\subset\mathbb{P}^{4}$, whose equation is
$\sum_{i=1}^{3}f_{i}g^{2}_{i}=0$ and which is smooth outside of a
curve $g_{1}=g_{2}=g_{3}=0$. Therefore, the 4-fold $X$ is not
rational due to
Proposition~\ref{proposition:fourfold-double-cover} and
Theorems~\ref{theorem:first} and \ref{theorem:six}.
\end{example}

It is natural to ask how many nodes can $X$ and $V$ have? The best
known upper bounds are due to \cite{Va83}. Namely,
$|\mathrm{Sing}(X)|\le\mathrm{A}_{3}(2r)$ and
$|\mathrm{Sing}(V)|\le\mathrm{A}_{4}(n)$, where
$\mathrm{A}_{i}(j)$ is a number of points
$(a_{1},\ldots,a_{i})\subset\mathbb{Z}^{i}$ such that
$$
(i-2){\frac{j}{2}}+1<\sum_{t=1}^{i}a_{t}\le {\frac{ij}{2}}
$$
and all $a_{t}\in(0,j)$. Hence, $|\mathrm{Sing}(X)|$ does not
exceed $68$, $180$ and $375$ when $r=3$, $4$ and $5$ respectively,
and $|\mathrm{Sing}(V)|$ does not exceed $45$, $135$ and $320$
when $n=4$, $5$ and $6$ respectively. In the case $n=4$ this bound
is sharp (see \cite{Fr86}). Moreover, there is only one nodal
quartic 3-fold with $45$ nodes (see \cite{JSV90}), so-called
Burkhardt quartic, which is rational and determinantal (see
\cite{Pet98}). In the case $r=3$ there is a better bound
$|\mathrm{Sing}(X)|\le 65$ which is sharp (see \cite{St78},
\cite{CaCe82}, \cite{Ba96}, \cite{JaRu97}, \cite{Wa98}). In the
case $n=5$ there are no known examples of nodal quntic
hypersurfaces in $\mathbb{P}^{4}$ with more than $130$ nodes (see
\cite{vSt93}).

\section{Preliminaries.}
\label{section:preliminaries}

Let $X$ be a variety and $B_{X}$ be a bounadry\footnote{Usually
boundaries are assumed to be effective (see \cite{KMM}), but we do
not assume this.}
 on $X$, i.e.
$B_X=\sum_{i=1}^{k}a_{i}B_{i}$, where $B_{i}$ is a prime divisor
on $X$ and $a_{i}\in\mathbb{Q}$. Basic notions, notations and
results related to the log pair $(X, B_{X})$ are contained in
\cite{KMM}, \cite{Ko91}, \cite{Ko97}. The log pair $(X, B_{X})$ is
called movable when every component $B_{i}$ is a linear system on
$X$ without fixed components. Basic properties of movable log
pairs are described in \cite{Al91}, \cite{Co95}, \cite{Co00},
\cite{Pu00}, \cite{Ch00a}, \cite{Ch03b}. In the following we
assume that $K_{X}$ and $B_{X}$ are $\mathbb{Q}$-Cartier divisors.

\begin{definition}
\label{definition:log-pull-back} A log pair $(V, B^{V})$ is called
a log pull back of the log pair $(X, B_{X})$ with respect to a
birational morphism $f:V\to X$ if
$K_{V}+B^{V}\sim_{\mathbb{Q}}f^{*}(K_{X}+B_{X})$ and
$$
B^{V}=f^{-1}(B_{X})-\sum_{i=1}^{n}a(X, B_{X}, E_{i})E_{i},
$$
where $E_{i}$ is an $f$-exceptional divisor and $a(X, B_{X},
E_{i})\in\mathbb{Q}$. The number $a(X, B_{X}, E_{i})$ is called a
discrepancy of the log pair $(X, B_X)$ in the $f$-exceptional
divisor $E_i$.
\end{definition}

\begin{definition}
\label{definition:log-resolution} A birational morphism $f:V\to X$
 is called a log resolution of the log pair $(X, B_{X})$ if the variety
$V$ is smooth and the union of all proper transforms of the
divisors $B_{i}$ and all $f$-exceptional divisors forms a divisor
with simple normal crossing.
\end{definition}

\begin{definition}
\label{difinition:center-of-log-canonical-singularities} A proper
irreducible subvariety $Y\subset X$ is called a center of log
canonical singularities of the log pair $(X, B_{X})$ if there are
a birational morphism $f:V\to X$ together with a not necessary
$f$-exceptional divisor $E\subset V$ such that $E$ is contained in
the support of the effective part of the divisor $\lfloor
B^{V}\rfloor$ and $f(E)=Y$. The set of all the centers of log
canonical sin\-gu\-la\-ri\-ties of the log pair $(X, B_{X})$ is
denoted by $\mathbb{LCS}(X, B_{X})$.
\end{definition}

\begin{definition}
\label{definition:log-canonical-singularities-subscheme} For a log
resolution $f:V\to X$ of $(X, B_{X})$ the subscheme
$\mathcal{L}(X, B_{X})$ associated to the ideal sheaf
$\mathcal{I}(X, B_{X})=f_{*}(\mathcal{O}_{V}(\lceil
-B^{V}\rceil))$ is called a log canonical singularity subscheme of
the log pair $(X, B_{X})$.
\end{definition}

The support of the log canonical singularity subscheme
$\mathcal{L}(X, B_{X})$ is a union of all elements in the set
$\mathbb{LCS}(X, B_{X})$. The following result is due to
\cite{Sh92} (see \cite{Ko91}, \cite{Am99}, \cite{Ch03b}).

\begin{theorem}
\label{theorem:Shokurov} Suppose that $B_{X}$ is effective and for
some nef and big divisor $H$ on $X$ the divisor $D=K_{X}+B_{X}+H$
is Cartier. Then $H^{i}(X, \mathcal{I}(X, B_{X})\otimes
\mathcal{O}_{X}(D))=0$ for $i>0$.
\end{theorem}

\begin{proof}
Let $f:W\to X$ be a log resolution of $(X,B_X)$. Then for $i>0$
$$
R^{i}f_{*}(f^{*}(\mathcal{O}_{X}(D))\otimes\mathcal{O}_{W}(\lceil-B^{W}\rceil))=0%
$$
by the relative Kawamata-Viehweg vanishing (see \cite{Ka82},
\cite{Vi82}). The equality
$$
R^{0}f_{*}(f^{*}(\mathcal{O}_{X}(D))\otimes\mathcal{O}_{W}(\lceil
-B^{W}\rceil))=\mathcal{I}(X, B_{X})\otimes \mathcal{O}_{X}(D)
$$
and the degeneration of local-to-global spectral sequence imply
that for all $i$
$$H^{i}(X, \mathcal{I}(X, B_{X})\otimes \mathcal{O}_{X}(D))=H^{i}(W, f^{*}(\mathcal{O}_{X}(D))\otimes\mathcal{O}_{W}(\lceil -B^{W}\rceil)),$$
but $H^{i}(W,
f^{*}(\mathcal{O}_{X}(D))\otimes\mathcal{O}_{W}(\lceil
-B^{W}\rceil))=0$ for $i>0$ by Kawamata-Viehweg vanishing.
\end{proof}

Consider the following application of
Theorem~\ref{theorem:Shokurov} (cf. \cite{EinLa93}, \cite{Ka97}).

\begin{lemma}
\label{lemma:non-vanishing} Let $\Sigma\subset\mathbb{P}^{n}$ be a
finite subset, $\mathcal{M}$ be a linear system of hypersurfaces
of degree $k$ passing through all points of the set $\Sigma$.
Suppose that the base locus of the linear system $\mathcal{M}$ is
zero-dimensional. Then the points of the set $\Sigma$ impose
independent linear conditions on the homogeneous forms on
$\mathbb{P}^{n}$ of degree $n(k-1)$.
\end{lemma}

\begin{proof}
Let $\Lambda\subset\mathbb{P}^{n}$ be a base locus of the linear
system $\mathcal{M}$. Then $\Sigma\subseteq\Lambda$ and $\Lambda$
is a finite subset in $\mathbb{P}^{n}$. Now consider sufficiently
general different divisors $H_{1},\ldots, H_{s}$ in the linear
system $\mathcal{M}$ for $s\gg 0$. Let $X=\mathbb{P}^{n}$ and
$B_{X}={\frac{n}{s}}\sum_{i=1}^{s}H_{i}$. Then
$\mathrm{Supp}(\mathcal{L}(X, B_{X}))=\Lambda$.

To prove the claim it is enough to prove that for every point
$P\in\Sigma$ there is a hypersurface in $\mathbb{P}^{n}$ of degree
$n(k-1)$ that passes through all the points in the set
$\Sigma\setminus P$ and does not pass through the point $P$. Let
$\Sigma\setminus P=\{P_{1},\ldots, P_{k}\}$, where $P_{i}$ is a
point of $X=\mathbb{P}^{n}$, and let $f:V\to X$ be a blow up at
the points of $\Sigma\setminus P$. Then
$$
K_{V}+(B_{V}+\sum_{i=1}^{k}(\mathrm{mult}_{P_{i}}(B_{X})-n)E_{i})+f^{*}(H)=f^{*}(n(k-1)H)-\sum_{i=1}^{k}E_{i},%
$$
where $E_{i}=f^{-1}(P_{i})$, $B_{V}=f^{-1}(B_{X})$ and $H$ is a
hyperplane in $\mathbb{P}^{n}$. By construction
$$
\mathrm{mult}_{P_{i}}(B_{X})=n\mathrm{mult}_{P_{i}}(\mathcal{M})\ge n%
$$
and the boundary
$\hat{B}_{V}=B_{V}+\sum_{i=1}^{k}(\mathrm{mult}_{P_{i}}(B_{X})-n)E_{i}$
is effective.

Let $\bar{P}=f^{-1}(P)$. Then $\bar{P}\in\mathbb{LCS}(V,
\hat{B}_{V})$ and $\bar{P}$ is an isolated center of log canonical
singularities of the log pair $(V, \hat{B}_{V})$, because in the
neighborhood of the point $P$ the birational morphism $f:V\to X$
is an isomorphism. On the other hand, the map
$$H^{0}(\mathcal{O}_{V}(f^{*}(n(k-1)H)-\sum_{i=1}^{k}E_{i}))\to H^{0}(\mathcal{O}_{\mathcal{L}(V, \hat{B}_{V})}\otimes\mathcal{O}_{V}(f^{*}(n(k-1)H)-\sum_{i=1}^{k}E_{i}))%
$$
is surjective by Theorem~\ref{theorem:Shokurov}. However, in the
neighborhood of the point $\bar{P}$ the support of the subscheme
$\mathcal{L}(V, \hat{B}_{V})$ consists just of the point
$\bar{P}$. The latter implies the existence of a divisor $D\in
|f^{*}(n(k-1)H)-\sum_{i=1}^{k}E_{i}|$ that does not pass through
$\bar{P}$. Thus, $f(D)$ is a hypersurface in $\mathbb{P}^{n}$ of
degree $n(k-1)$ that passes through the points of $\Sigma\setminus
P$ and does not pass through the point $P\in\Sigma$.
\end{proof}

Actually, the proof of Lemma~\ref{lemma:non-vanishing} implies
Theorem~\ref{theorem:six}.

\begin{proof}[Proof of Theorem~\ref{theorem:six}]
We have a double cover $\pi:X\to\mathbb{P}^{3}$ branched over a
nodal hypersurface $S\subset\mathbb{P}^{3}$ of degree $2r$, a
linear subsystem
$\mathcal{H}\subset|\mathcal{O}_{\mathbb{P}^{3}}(k)|$ of
hypersurfaces vanishing at the points of the set
$\mathrm{Sing}(S)$ for $k<r$, a nodal hypersurface
$V\subset\mathbb{P}^{4}$ of degree $n$, a linear subsystem
$\mathcal{D}\subset|\mathcal{O}_{\mathbb{P}^{4}}(l)|$ of
hypersurfaces vanishing at $\mathrm{Sing}(V)$ for
$l<{\frac{n}{2}}$, and
$$
\mathrm{dim}(\mathrm{Bs}(\hat{\mathcal{H}}))=\mathrm{dim}(\mathrm{Bs}(\hat{\mathcal{D}}))=0,
$$
where $\hat{\mathcal{H}}=\mathcal{H}\vert_{S}$ and
$\hat{\mathcal{D}}=\mathcal{D}\vert_{V}$. Due to
Proposition~\ref{proposition:defect} we must show that the nodes
of the surface $S\subset\mathbb{P}^{3}$ and the nodes of the
hypersurface  $X\subset\mathbb{P}^{4}$ impose independent linear
conditions on homogeneous forms of degree $3r-4$ and $2n-5$
respectively.

Suppose that the stronger condition
$$
\mathrm{dim}(\mathrm{Bs}({\mathcal{H}}))=\mathrm{dim}(\mathrm{Bs}({\mathcal{D}}))=0
$$
holds, which is enough for Corollary~\ref{corollary:theorem-six}.
Then Lemma~\ref{lemma:non-vanishing} immediately implies that the
nodes of $S$ and the nodes of $X$ impose independent linear
conditions on homogeneous forms of degree $3r-4$ and $2n-5$
respectively. In the general case we can repeat the proof of
Lemma~\ref{lemma:non-vanishing} interchanging the boundary
${\frac{3}{s}}\sum_{i=1}^{s}H_{i}$ with the boundary
$S+{\frac{1}{s}}\sum_{i=1}^{s}H_{i}$ for the proof of the
$\mathbb{Q}$-factoriality of $X$ and interchanging the boundary
${\frac{4}{s}}\sum_{i=1}^{s}H_{i}$ with the boundary
$V+{\frac{2}{s}}\sum_{i=1}^{s}H_{i}$ for the proof of the
$\mathbb{Q}$-factoriality of the 3-fold $V$.
\end{proof}

The following result is Theorem~17.4 in \cite{Ko91} and
Theorem~7.4 in \cite{Ko97}.

\begin{theorem}
\label{theorem:connectedness} Let $g:X\to Z$ be a contraction such
that $g_{*}(\mathcal{O}_{X})=\mathcal{O}_{Z}$, $h:V\to X$ be a log
resolution of the log pair $(X, B_{X})$. Suppose the divisor
$-(K_{X}+B_{X})$ is $g$-nef and $g$-big, and
$\mathrm{codim}(g(B_{i})\subset Z)\geq 2$ whenever $a_{i}<0$.
Define $a_{E}\in\mathbb{Q}$ by means of the equivalence
$$
K_{V}\sim_\mathbb{Q} f^{*}(K_{X}+B_{X})+\sum_{E\subset V} a_{E}E,
$$
where $E\subset V$ is a not necessary $h$-exceptional divisor.
Then the locus $\cup_{a_{E}\le -1}E$ is connected in a
neighborhood of every fiber of the morphism $g\circ h$.
\end{theorem}

The following result is a corollary of
Theorem~\ref{theorem:connectedness} (see  Theorem~17.6 in
\cite{Ko91}).

\begin{theorem}
\label{theorem:log-adjunction} Suppose that $B_{X}$ is effective
and $\lfloor B_{X}\rfloor=\emptyset$. Let $S\subset X$ be an
effective irreducible divisor such that the divisor
$K_{X}+S+B_{X}$ is $\mathbb{Q}$-Cartier. Then $(X, S+B_{X})$ is
purely log terminal if and only if $(S, \mathrm{Diff}_{S}(B_{X}))$
is Kawamata log terminal.
\end{theorem}

\begin{definition}
\label{definition:center-of-canonical-singularities} A proper
irreducible subvariety $Y\subset X$ is called a center of
canonical singularities of $(X, B_{X})$ if there is a birational
morphism $f:W\to X$ and an $f$-ex\-cep\-tional divisor $E\subset
W$ such that the discrepancy $a(X, B_{X}, E)\leq 0$ and $f(E)=Y$.
The set of all centers of canonical sin\-gu\-la\-ri\-ties of the
log pair $(X, B_{X})$ is denoted by $\mathbb{CS}(X, B_{X})$.
\end{definition}

\begin{corollary}
\label{corollary:log-adjunction-I} Let $H$ be an effective Cartier
divisor on $X$ and $Z\in \mathbb{CS}(X, B_{X})$, suppose that $X$
and $H$ are smooth in the generic point of $Z$, $Z\subset H$,
$H\not\subset\mathrm{Supp}(B_{X})$ and $B_{X}$ is an effective
boundary. Then $\mathbb{LCS}(H, B_{X}\vert_{H})\ne\emptyset$.
\end{corollary}

The following result is Corollary 7.3 in \cite{Pu00}, which holds
even over fields of positive characteristic and implicitly goes
back to \cite{IsMa71} (see \cite{Co95}, \cite{Co00},
\cite{Kaw01}).

\begin{theorem}
\label{theorem:Iskovskikh} Suppose that $X$ is smooth,
$\mathrm{dim}(X)\ge 3$, the boundary $B_{X}$ is effective and
movable, and the set $\mathbb{CS}(X, M_{X})$ contains a closed
point $O\in X$. Then $\mathrm{mult}_{O}(B_{X}^{2})\geq 4$ and the
equality implies $\mathrm{mult}_{O}(B_{X})=2$ and
$\mathrm{dim}(X)=3$.
\end{theorem}

The following result is implied by Theorem~3.10 in \cite{Co00} and
Theorem~\ref{theorem:log-adjunction}.

\begin{theorem}
\label{theorem:Corti} Suppose that $\mathrm{dim}(X)\ge 3$, $B_{X}$
is effective, and the set $\mathbb{CS}(X, B_{X})$ contains an
ordinary double point $O$ of $X$. Then
$\mathrm{mult}_{O}(B_{X})\geq 1$, where $\mathrm{mult}_{O}(B_{X})$
is defined by means of the regular blow up of $O$. Moreover,
$\mathrm{mult}_{O}(B_{X})=1$ implies $\mathrm{dim}(X)=3$.
\end{theorem}

The following result is an easy modification of
Theorem~\ref{theorem:Corti}.

\begin{proposition}
\label{proposition:non-simple-double-point} Suppose that
$\mathrm{dim}(X)=3$, $B_{X}$ is effective, and the set
$\mathbb{CS}(X, B_{X})$ contains an isolated singular point $O$ of
the variety $X$, which is locally isomorphic to the singularity
$y^{3}=\sum_{i=1}^{3}x_{i}^{2}$. Let $f:W\to X$ be a blow up of
$O$ and $\mathrm{mult}_{O}(B_{X})$ be a rational number defined by
means of the equivalence
$$
f^{-1}(B_{X})\sim_{\mathbb{Q}}f^{-1}(B_{X})+\mathrm{mult}_{O}(B_{X})E,
$$
where $E$ is an $f$-exceptional divisor. Then
$\mathrm{mult}_{O}(B_{X})\ge {\frac{1}{2}}$.
\end{proposition}

\begin{proof} The 3-fold $W$ is
smooth, $E$ is isomorphic to a cone in $\mathbb{P}^{3}$ over a
smooth conic, the restriction $-E\vert_{E}$ is rationally
equivalent to a hyperplane section of $E\subset\mathbb{P}^{3}$,
and
$$
K_{W}+B_{W}\sim_{\mathbb{Q}}
f^{*}(K_{X}+B_{X})+(1-\mathrm{mult}_{O}(B_{X}))E,
$$
where $B_{W}=f^{-1}(B_{X})$. Suppose that
$\mathrm{mult}_{O}(B_{X})<{\frac{1}{2}}$. Then there is a proper
irreducible subvariety $Z\subset E$ such that $Z\in\mathbb{CS}(W,
B_{W})$. Hence, $\mathbb{LCS}(E, B_{W}\vert_{E})\ne\emptyset$ by
Theorem~\ref{theorem:log-adjunction}.

Let $B_{E}=B_{W}\vert_{E}$. Then $\mathbb{LCS}(E, B_{E})$ does not
contains curves on $E$, because otherwise the intersection of
$B_{E}$ with the ruling of $E$ is greater than ${\frac{1}{2}}$,
which is impossible due to our assumption
$\mathrm{mult}_{O}(B_{X})<{\frac{1}{2}}$. Therefore,
$\mathrm{dim}(\mathrm{Supp}(\mathcal{L}(E, B_{E})))=0$.

Let $H$ be a hyperplane of $E\subset\mathbb{P}^{3}$. Then
$$
K_{E}+B_{E}+(1-\mathrm{mult}_{O}(B_{X}))H\sim_{\mathbb{Q}}-H
$$
and $H^{0}(\mathcal{O}_{E}(-H))=0$. On the other hand,
Theorem~\ref{theorem:Shokurov} implies surjectivity
$$
H^{0}(\mathcal{O}_{E}(-H))\to H^{0}(\mathcal{O}_{\mathcal{L}(E, B_{E})})\to 0,%
$$
which is a contradiction.
\end{proof}

The following result is due to \cite{Co95} (see \cite{Pu00},
\cite{Ch03b}).

\begin{theorem}
\label{theorem:Nother-Fano} Let $X$ be a Fano variety with
$\mathrm{Pic}(X)\cong \mathbb{Z}$ with terminal
$\mathbb{Q}$-factorial singularities such that either $X$ is not
birationally rigid or $\mathrm{Bir}(X)\ne\mathrm{Aut}(X)$. Then
there are a linear system $\mathcal{M}$ on $X$ having no fixed
components and $\mu\in\mathbb{Q}_{>0}$ such that the singularities
of the movable log pair $(X, \mu\mathcal{M})$ are not canonical
and $\mu\mathcal{M}\sim_{\mathbb{Q}} -K_{X}$.
\end{theorem}

The following result is proved in \cite{Bes83} using the vanishing
theorem for $2$-connected divisors on algebraic surfaces in
\cite{Ram72} in a way similar to \cite{Bom73} and \cite{vdVe79}.

\begin{theorem}
\label{theorem:Bese} Let $\pi:Y\to\mathbb{P}^2$ be the blow up at
points $P_{1},\ldots, P_{s}$ on $\mathbb{P}^2$,
$s\leq{\frac{d^{2}+9d+10}{6}}$, at most $k(d+3-k)-2$ of the points
$P_i$ lie on a curve of degree $k\leq {\frac{d+3}{2}}$, where
$d\ge 3$ is a natural number. Then the linear system
$|\pi^{*}(\mathcal{O}_{\mathbb{P}^2}(d))-\sum_{i=1}^{s}E_{i}|$ is
free, where $E_{i}=\pi^{-1}(P_{i})$
\end{theorem}

\begin{corollary}
\label{corollary:Bese} Let $\Sigma\subset\mathbb{P}^{2}$ be a
finite subset such that the inequality $|\Sigma|\leq
{\frac{d^{2}+9d+16}{6}}$ holds and at most $k(d+3-k)-2$ points of
$\Sigma$ lie on a curve of degree $k\leq {\frac{d+3}{2}}$, where
$d\ge 3$ is a natural number. Then for every point $P\in\Sigma$
there is a curve $C\subset\mathbb{P}^{2}$ of degree $d$ that
passes through all the points in $\Sigma\setminus P$ and does not
pass through the point $P$.
\end{corollary}

In the case $d=3$ the claim of Theorem~\ref{theorem:Bese} is
nothing but the freeness of the anticanonical linear system of a
weak del Pezzo surface of degree $9-s\ge 2$ (see \cite{De80}).

\section{Double solids.}
\label{section:double-solids}

In this section we prove Theorem~\ref{theorem:first}. Let
$\pi:X\to\mathbb{P}^{3}$ be a double cover  branched over a nodal
hypersurface $S\subset\mathbb{P}^{3}$ of degree $2r$ such that
$|\mathrm{Sing}(S)|\le{\frac{(2r-1)r}{3}}$. We must show that the
nodes of $S\subset\mathbb{P}^{3}$ impose independent linear
conditions on homogeneous forms of degree $3r-4$ on
$\mathbb{P}^{3}$ due to Proposition~\ref{proposition:defect} and
Remark~\ref{remark:factoriality}. Moreover, we may assume $r\ge
3$, because in the case $r\le 2$ the required claim is trivial.

\begin{definition}
\label{definition:general-position} The points of a subset
$\Gamma\subset{\mathbb P}^{s}$ satisfy the property $\nabla$ if at
most $k(2r-1)$ points of the set $\Gamma$ can lie on a curve in
$\mathbb{P}^{s}$ of degree $k\in\mathbb{N}$.
\end{definition}

Let $\Sigma=\mathrm{Sing}(S)\subset\mathbb{P}^{3}$.

\begin{proposition}
\label{proposition:nodes-in-general-position} The points of the
subset $\Sigma\subset\mathbb{P}^{3}$ satisfy the property
$\nabla$.
\end{proposition}

\begin{proof}
Let $L\subset{\mathbb P}^3$ be a line and $\Pi\subset{\mathbb
P}^3$ be a sufficiently general hyperplane passing through the
line $L$. Then $\Pi\not\subset S$ and $\Pi\cap S=L\cup Z$, where
$Z\subset\Pi$ is a plane curve of degree $2r-1$. Moreover, we have
$$
\Sigma\cap L=\mathrm{Sing}(S)\cap L\subseteq L\cap Z,
$$
but $|L\cap Z|\le 2r-1$. Thus, at most $2r-1$ points of $\Sigma$
can lie on a line.

Let $C\subset{\mathbb P}^3$ be a curve of degree $k>1$. We must
show that at most $k(2r-1)$ points of $\Sigma$ can lie on $C$. We
may assume that $C$ is irreducible and reduced. Consider a general
cone $Y\subset{\mathbb P}^3$ over the curve $C$. Then
$Y\not\subset S$ and $Y\cap S=C\cup R$, where $R$ is an
irreducible reduced curve of degree $k(2r-1)$. As above we have
the inclusion
$$
\Sigma\cap C=\mathrm{Sing}(S)\cap C\subseteq C\cap R,
$$
but in the set-theoretic sense $|C\cap R|\le (2r-1)k$. Hence, at
most $k(2r-1)$ points of the subset $\Sigma\subset\mathbb{P}^{3}$
can lie on the irreducible reduced curve $C\subset\mathbb{P}^{3}$
of degree $k$.
\end{proof}

Fix a point $P\in\Sigma$. To prove that the points of
$\Sigma\subset\mathbb{P}^{3}$ impose independent linear conditions
on homogeneous forms of degree $3r-4$ it is enough to construct a
hypersurface in $\mathbb{P}^{3}$ of degree $3r-4$ that passes
through $\Sigma\setminus P$ and does not pass through
$P\in\Sigma$.

\begin{lemma}
\label{lemma:double-solid-I} Suppose $\Sigma\subset\Pi$ for some
hyperplane $\Pi\subset\mathbb{P}^{3}$. Then there is a
hypersurface in $\mathbb{P}^{3}$ of degree $3r-4$ that passes
through $\Sigma\setminus P$ and does not pass through
$P\in\Sigma$.
\end{lemma}

\begin{proof}
Let us apply Corollary~\ref{corollary:Bese} to $\Sigma\subset\Pi$
and $d=3r-4\ge 5$. We must check that all the conditions of
Corollary~\ref{corollary:Bese} are satisfied, which is easy but
not obvious. First of all
$$
|\Sigma|\le {\frac{(2r-1)r}{3}}\Rightarrow |\Sigma|\leq {\frac{d^{2}+9d+16}{6}}%
$$
and at most $d=3r-4$ points of $\Sigma$ can lie on a line in $\Pi$
because $r\ge 3$ and the points of the subset $\Sigma\subset\Pi$
satisfy the property $\nabla$ due to
Proposition~\ref{proposition:nodes-in-general-position}.

Now we must prove that at most $k(3r-1-k)-2$ points of $\Sigma$
can lie on a curve of degree $k\leq {\frac{3r-1}{2}}$. The case
$k=1$ is already done. Moreover, at most $k(2r-1)$ points of the
set $\Sigma$ can lie on a curve of degree $k$ by
Proposition~\ref{proposition:nodes-in-general-position}. Thus, we
must show that
$$
k(3r-1-k)-2\ge k(2r-1)
$$
for all $k\le {\frac{3r-1}{2}}$. Moreover, we must prove the
latter inequality only for such $k>1$ that the inequality
$k(3r-1-k)-2<|\Sigma|$ holds, because otherwise the corresponding
condition on the points of the set $\Sigma$ is vacuous. Moreover,
we have
$$
k(3r-1-k)-2\ge k(2r-1)\iff r>k,
$$
because $k>1$. Suppose that the inequality $r\le k$ holds for some
natural number $k$ such that $k\le {\frac{3r-1}{2}}$ and
$k(3r-1-k)-2<|\Sigma|$. Let $g(x)=x(3r-1-x)-2$. Then $g(x)$ is
increasing for $x<{\frac{3r-1}{2}}$. Thus, we have $g(k)\ge g(r)$,
because ${\frac{3r-1}{2}}\ge k\ge r$. Hence,
$$
{\frac{(2r-1)r}{3}}\ge|\Sigma|>g(k)\ge g(r)=r(2r-1)-2,
$$
which is impossible when $r\ge 3$.

Therefore, there is a curve $C\subset\Pi$ of degree $3r-4$ that
passes trough $\Sigma\setminus P$ and does not pass through $P$ by
Corollary~\ref{corollary:Bese}. Let $Y\subset\mathbb{P}^{3}$ be a
sufficiently general cone over the curve
$C\subset\Pi\cong\mathbb{P}^{2}$. Then $Y\subset\mathbb{P}^{3}$ is
a hypersurface of degree $3r-4$ that passes through all the points
of the set $\Sigma\setminus P$ and does not pass through the point
$P\in\Sigma$.
\end{proof}

Take a sufficiently general hyperplane $\Pi\subset\mathbb{P}^{3}$.
Let $\psi:\mathbb{P}^{3}\dasharrow\Pi$ be a projection from a
sufficiently general point $O\in\mathbb{P}^{3}$,
$\Sigma^{\prime}=\psi(\Sigma)\subset\Pi\cong\mathbb{P}^{2}$ and
$\hat{P}=\psi(P)\in\Sigma^{\prime}$.

\begin{lemma}
\label{lemma:double-solid-II} Suppose that the points of
$\Sigma^{\prime}\subset\Pi$ satisfy the property $\nabla$. Then
there is a hypersurface in $\mathbb{P}^{3}$ of degree $3r-4$
containing $\Sigma\setminus P$ and not passing through $P$.
\end{lemma}

\begin{proof}
The proof of the claim of Lemma~\ref{lemma:double-solid-I} implies
the existence of a curve $C\subset\Pi$ of degree $3r-4$ that
passes through $\Sigma^{\prime}\setminus\hat{P}$ and does not pass
through $\hat{P}$. Let $Y\subset\mathbb{P}^{3}$ be a cone over the
curve $C$ with the vertex $O$. Then $Y\subset\mathbb{P}^{3}$ is a
hypersurface of degree $3r-4$ that passes through $\Sigma\setminus
P$ and does not pass through the point $P\in\Sigma$.
\end{proof}

It seems to us that the points of the subset
$\Sigma^{\prime}\subset\Pi$ always satisfy the property $\nabla$
due to the generality in the choice of the projection
$\psi:\mathbb{P}^{3}\dasharrow\Pi$. Unfortunately, we are unable
to prove it. Hence, we may assume that the points of the subset
$\Sigma^{\prime}\subset\Pi\cong\mathbb{P}^{2}$ do not satisfy the
property $\nabla$. Let us clarify this assumption.

\begin{definition}
\label{definition:general-position-restricted} The points of a
subset $\Gamma\subset{\mathbb P}^{s}$ satisfy the property
$\nabla_{k}$ if at most $i(2r-1)$ points of the set $\Gamma$ can
lie on a curve in $\mathbb{P}^{s}$ of degree $i\in\mathbb{N}$ for
all $i\le k$.
\end{definition}

Therefore, there is a smallest $k\in\mathbb{N}$ such that the
points of $\Sigma^{\prime}\subset\Pi$ do not satisfy the property
$\nabla_{k}$. Namely, there is a subset
$\Lambda_{k}^{1}\subset\Sigma$ such that
$|\Lambda_{k}^{1}|>k(2r-1)$ and all the points of the set
$$
\tilde{\Lambda}_{k}^{1}=\psi(\Lambda_{k}^{1})\subset\Sigma^{\prime}\subset\Pi\cong\mathbb{P}^{2}
$$
lie on a curve $C\subset\Pi$ of degree $k$. Moreover, the curve
$C$ is irreducible and reduced due to the minimality of $k$. In
the case when the points of the subset
$\Sigma^{\prime}\setminus\tilde{\Lambda}_{k}^{1}\subset\Pi$ does
not satisfy the property $\nabla_{k}$ we can find subset
$\Lambda_{k}^{2}\subset\Sigma\setminus\Lambda_{k}^{1}$ such that
$|\Lambda_{k}^{2}|>k(2r-1)$ and all the points of the set
$\tilde{\Lambda}_{k}^{2}=\psi(\Lambda_{k}^{2})$ lie on an
irreducible curve of degree $k$. Thus, we can iterate this
construction $c_{k}$ times and get $c_{k}>0$ disjoint subsets
$$
\Lambda_{k}^{i}\subset\Sigma\setminus\bigcup_{j=1}^{i-1}\Lambda_{k}^{j}\subsetneq\Sigma
$$
such that $|\Lambda_{k}^{i}|>k(2r-1)$, all the points of the
subset
$\tilde{\Lambda}_{k}^{i}=\psi(\Lambda_{k}^{i})\subset\Sigma^{\prime}$
lie on an irreducible reduced curve on $\Pi$ of degree $k$, and
all the points of the subset
$$
\Sigma^{\prime}\setminus\bigcup_{i=1}^{c_{k}}\tilde{\Lambda}_{k}^{i}\subset\Pi\cong\mathbb{P}^{2}
$$
satisfy the property $\nabla_{k}$. Now we can repeat this
construction for the property $\nabla_{k+1}$ and find $c_{k+1}\ge
0$ disjoint subsets
$$
\Lambda_{k+1}^{i}\subset(\Sigma\setminus\bigcup_{i=1}^{c_{k}}\Lambda_{k}^{i})\setminus\bigcup_{j=1}^{i-1}\Lambda_{k+1}^{j}\subset\Sigma\setminus\bigcup_{i=1}^{c_{k}}\Lambda_{k}^{i}\subsetneq\Sigma
$$
such that $|\Lambda_{k+1}^{i}|>(k+1)(2r-1)$, the points of
$\tilde{\Lambda}_{k+1}^{i}=\psi(\Lambda_{k+1}^{i})\subset\Sigma^{\prime}$
lie on an irreducible reduced curve on $\Pi$ of degree $k+1$, and
the points of the subset
$$
\Sigma^{\prime}\setminus\bigcup_{j=k}^{k+1}\bigcup_{i=1}^{c_{k}}\tilde{\Lambda}_{j}^{i}\subsetneq\Sigma^{\prime}\subset\Pi\cong\mathbb{P}^{2}
$$
satisfy the property $\nabla_{k+1}$. Now we can iterate this
construction for $\nabla_{k+2}, \ldots, \nabla_{l}$ and get
disjoint subsets $\Lambda_{j}^{i}\subset\Sigma$ for
$j=k,\ldots,l\ge k$ such that $|\Lambda_{j}^{i}|>j(2r-1)$, all the
points of the subset
$\tilde{\Lambda}_{j}^{i}=\psi(\Lambda_{j}^{i})\subset\Sigma^{\prime}$
lie on an irreducible reduced curve of degree $j$ in $\Pi$, and
all the points in the subset
$$
\bar{\Sigma}=\Sigma^{\prime}\setminus\bigcup_{j=k}^{l}\bigcup_{i=1}^{c_{j}}\tilde{\Lambda}_{j}^{i}\subsetneq\Sigma^{\prime}\subset\Pi\cong\mathbb{P}^{2}
$$
satisfy the property $\nabla$, where $c_{j}\ge 0$ is a number of
subsets $\tilde{\Lambda}_{j}^{i}$. Note, that $c_{k}>0$.

\begin{remark}
\label{remark:emptyness-of-subsets} The subset
$\Lambda_{k}^{1}\subset\Sigma$ is non-empty. However, every subset
$\Lambda_{j}^{i}\subset\Sigma$ a priori can be empty when $j\ne k$
or $i\ne 1$. Moreover, the subset
$\bar{\Sigma}\subset\Sigma^{\prime}$ can be empty as well.
\end{remark}

\begin{remark}
\label{remark:number-of-good-points} The inequality
$|\bar{\Sigma}|<{\frac{(2r-1)r}{3}}-\sum_{i=k}^{l}c_{i}(2r-1)i={\frac{(2r-1)}{3}}(r-3\sum_{i=k}^{l}ic_{i})$
holds.
\end{remark}

\begin{corollary}
\label{corollary:from-the-number-of-good-points} The inequality
$\sum_{i=k}^{l}ic_{i}<{\frac{r}{3}}$ holds.
\end{corollary}

In particular, $\Lambda_{j}^{i}\ne\emptyset$ implies
$j<{\frac{r}{3}}$.

\begin{lemma}
\label{lemma:double-solid-III} Suppose that
$\Lambda_{j}^{i}\ne\emptyset$. Let $\mathcal{M}$ be a linear
system of hypersurfaces of degree $j$ in $\mathbb{P}^{3}$ passing
through all the points in $\Lambda_{j}^{i}$. Then the base locus
of $\mathcal{M}$ is zero-dimensional.
\end{lemma}

\begin{proof}
By the construction of the set $\Lambda_{j}^{i}$ all the points of
the subset
$$
\tilde{\Lambda}_{j}^{i}=\psi(\Lambda_{j}^{i})\subset\Sigma^{\prime}\subset\Pi\cong\mathbb{P}^{2}
$$
lie on an irreducible reduced curve $C\subset\Pi$ of degree $j$.
Let $Y\subset\mathbb{P}^{3}$ be a cone over $C$ with the vertex
$O$. Then $Y$ is a hypersurface in $\mathbb{P}^{3}$ of degree $j$
that contains all the points of the set $\Lambda_{j}^{i}$.
Therefore, $Y\in\mathcal{M}$.

Suppose that the base locus of the linear system $\mathcal{M}$
contains an irreducible reduced curve $Z\subset\mathbb{P}^{3}$.
Then $Z\subset Y$ and $\psi(Z)=C$. Moreover,
$\Lambda_{j}^{i}\subset Z$, because $\Lambda_{j}^{i}\not\subset Z$
implies that $\tilde{\Lambda}_{j}^{i}\not\subset C$ due to the
generality of $\psi$. Finally, the restriction $\psi\vert_{Z}:Z\to
C$ is a birational morphism, because the projection $\psi$ is
general. Hence, $\mathrm{deg}(Z)=j$ and  $Z$ contains at least
$|\Lambda_{j}^{i}|>j(2r-1)$ points of $\Sigma$. The latter
contradicts
Proposition~\ref{proposition:nodes-in-general-position}.
\end{proof}

\begin{corollary}
\label{corollary:k-is-at-most-two} The inequality $k\ge 2$ holds.
\end{corollary}

For every $\Lambda_{j}^{i}\ne\emptyset$ let
$\Xi_{j}^{i}\subset\mathbb{P}^{3}$ be a base locus of the linear
system of hypersurfaces of degree $j$ in $\mathbb{P}^{3}$ passing
through all the points in $\Lambda_{j}^{i}$. For
$\Lambda_{j}^{i}=\emptyset$ put $\Xi_{j}^{i}=\emptyset$. Then
$\Xi_{j}^{i}$ is a finite set by
Lemma~\ref{lemma:double-solid-III} and
$\Lambda_{j}^{i}\subseteq\Xi_{j}^{i}$ by construction.

\begin{lemma}
\label{lemma:double-solid-V} Suppose that
$\Xi_{j}^{i}\ne\emptyset$. Then the points of the subset
$\Xi_{j}^{i}\subset\mathbb{P}^{3}$ impose independent linear
conditions on homogeneous forms on $\mathbb{P}^{3}$ of degree
$3(j-1)$.
\end{lemma}

\begin{proof}
The claim follows from Lemma~\ref{lemma:non-vanishing}.
\end{proof}

\begin{corollary}
\label{corollary:bad-good-points} Suppose that
$\Lambda_{j}^{i}\ne\emptyset$. Then the points of the subset
$\Lambda_{j}^{i}\subset\mathbb{P}^{3}$ impose independent linear
conditions on homogeneous forms on $\mathbb{P}^{3}$ of degree
$3(j-1)$.
\end{corollary}

\begin{lemma}
\label{lemma:double-solid-VI} Suppose that
$\bar{\Sigma}=\emptyset$. Then there is a hypersurface in
$\mathbb{P}^{3}$ of degree $3r-4$ containing $\Sigma\setminus P$
and not passing through the point $P$.
\end{lemma}

\begin{proof}
The set $\Sigma$ is a disjoint union
$\cup_{j=k}^{l}\cup_{i=1}^{c_{j}}\Lambda_{j}^{i}$ and there is a
unique set $\Lambda_{a}^{b}$ containing the point $P$. In
particular, $P\in\Xi_{a}^{b}$. On the other hand, the union
$\cup_{j=k}^{l}\cup_{i=1}^{c_{j}}\Xi_{j}^{i}$ is not necessary
disjoint. Thus, a priori the point $P$ can be contained in many
sets $\Xi_{j}^{i}$.

For every $\Xi_{j}^{i}\ne\emptyset$ containing $P$ there is a
hypersurface of degree $3(j-1)$ that passes through
$\Xi_{j}^{i}\setminus P$ and does not pass through  $P$ by
Lemma~\ref{lemma:double-solid-V}. For every
$\Xi_{j}^{i}\ne\emptyset$ not containing the point $P$ there is a
hypersurface of degree $j$ that passes through $\Xi_{j}^{i}$ and
does not pass through the point $P$ by the definition of the set
$\Xi_{j}^{i}$. Moreover, $j<3(j-1)$, because $k\ge 2$ by
Corollary~\ref{corollary:k-is-at-most-two}. Therefore, for every
$\Xi_{j}^{i}\ne\emptyset$ there is a hypersurface
$F_{i}^{j}\subset\mathbb{P}^{3}$ of degree $3(j-1)$ that passes
through $\Xi_{j}^{i}\setminus P$ and does not pass through the
point $P$. Let
$$
F=\bigcup_{j=k}^{l}\bigcup_{i=1}^{c_{j}}F_{j}^{i}\subset\mathbb{P}^{3}
$$
be a possibly reducible hypersurface of degree
$\sum_{i=k}^{l}3(i-1)c_{i}$. Then $F$ passes through all the
points of the set $\Sigma\setminus P$ and does not pass through
the point $P$. Moreover, we have
$$
\mathrm{deg}(F)=\sum_{i=k}^{l}3(i-1)c_{i}<\sum_{i=k}^{l}3ic_{i}<r<3r-4
$$
by Corollary~\ref{corollary:from-the-number-of-good-points}, which
implies the claim.
\end{proof}

Let $\hat{\Sigma}=\cup_{j=k}^{l}\cup_{i=1}^{c_{j}}\Lambda_{j}^{i}$
and $\check{\Sigma}=\Sigma\setminus\hat{\Sigma}$. Then
$\Sigma=\hat{\Sigma}\cup\check{\Sigma}$ and
$\psi(\check{\Sigma})=\bar{\Sigma}\subset\Pi$. Therefore, we
proved Theorem~\ref{theorem:first} in the extreme cases:
$\hat{\Sigma}=\emptyset$ and $\check{\Sigma}=\emptyset$. Now we
must combine the proofs of the Lemmas~\ref{lemma:double-solid-II}
and \ref{lemma:double-solid-VI} to prove
Theorem~\ref{theorem:first} in the case $\hat{\Sigma}\ne\emptyset$
and $\check{\Sigma}\ne\emptyset$.

\begin{remark}
\label{remark:hypersurface-through-bad-points} The proof of
Lemma~\ref{lemma:double-solid-VI} implies the existence of a
hypersurface $F\subset\mathbb{P}^{3}$ of degree
$\sum_{i=k}^{l}3(i-1)c_{i}$ that passes through all the points of
the subset $\hat{\Sigma}\setminus P\subsetneq\Sigma$ and does not
pass through the point $P\in\Sigma$.
\end{remark}

Put $d=3r-4-\sum_{i=k}^{l}3(i-1)c_{i}$. Let us check that the
subset $\bar{\Sigma}\subset\Pi\cong\mathbb{P}^{2}$ and the number
$d$ satisfy all the conditions of Theorem~\ref{theorem:Bese}. We
may assume that $\emptyset\ne\hat{\Sigma}\subsetneq\Sigma$.

\begin{lemma}
\label{lemma:double-solid-VII} The inequality $d\ge 6$ holds.
\end{lemma}

\begin{proof}
The inequality
$$
|\bar{\Sigma}|<{\frac{(2r-1)r}{3}}-\sum_{i=k}^{l}c_{i}(2r-1)i={\frac{(2r-1)}{3}}(r-3\sum_{i=k}^{l}ic_{i})
$$
implies $\sum_{i=k}^{l}3ic_{i}<r$. Thus,
$d=3r-4-\sum_{i=k}^{l}3(i-1)c_{i}>2r-4+3c_{k}\ge 2r-1\ge 5$.
\end{proof}

\begin{lemma}
\label{lemma:double-solid-VIII} The inequality $|\bar{\Sigma}|\leq
{\frac{d^{2}+9d+10}{6}}$ holds.
\end{lemma}

\begin{proof}
By construction
$|\bar{\Sigma}|<{\frac{(2r-1)}{3}}(r-3\sum_{i=k}^{l}ic_{i})$.
Thus, we must show that
$$
2(2r-1)(r-3\sum_{i=k}^{l}ic_{i})\le
(3r-4-\sum_{i=k}^{l}3(i-1)c_{i})^{2}+9(3r-4-\sum_{i=k}^{l}3(i-1)c_{i})+10,
$$
where $c_{k}\ge 1$ and $\sum_{i=k}^{l}3ic_{i}<r$. However, we have
$$
(3r-4-\sum_{i=k}^{l}3(i-1)c_{i})^{2}+9(3r-4-\sum_{i=k}^{l}3(i-1)c_{i})+10>
(2r-4+3c_{k})^{2}+9(2r-4+3c_{k})+10
$$
and
$$
(2r-4+3c_{k})^{2}+9(2r-4+3c_{k})+10\ge
(2r-1)^{2}+9(2r-1)+10=4r^{2}+14r+2,
$$
which implies
$4r^{2}+14r+2>4r^{2}-2r>2(2r-1)(r-3\sum_{i=k}^{l}ic_{i})$.
\end{proof}

\begin{lemma}
\label{lemma:double-solid-IX} At most $d$ of the points of the
subset $\bar{\Sigma}\subset\mathbb{P}^{2}$ lie on a line in
$\mathbb{P}^{2}$.
\end{lemma}

\begin{proof}
The points of the subset $\bar{\Sigma}\subset\mathbb{P}^{2}$
satisfy the property $\nabla$. In particular, no more than $2r-1$
of the points of $\bar{\Sigma}$ lie on a line in $\mathbb{P}^{2}$.
On the other hand, we have
$$
d=3r-4-\sum_{i=k}^{l}3(i-1)c_{i}>2r-4+3c_{k}\ge 2r-1,
$$
which implies the claim.
\end{proof}

\begin{lemma}
\label{lemma:double-solid-X} At most $k(d+3-k)-2$ points of
$\bar{\Sigma}$ lie on a curve in $\mathbb{P}^{2}$ of degree $k\leq
{\frac{d+3}{2}}$.
\end{lemma}

\begin{proof}
We may assume that $k>1$ due to Lemma~\ref{lemma:double-solid-IX}.
The points of the subset $\bar{\Sigma}\subset\mathbb{P}^{2}$
satisfy the property $\nabla$. Thus, at most $(2r-1)k$ of the
points of $\bar{\Sigma}$ lie on a curve in $\mathbb{P}^{2}$ of
degree $k$. Therefore, to conclude the proof it is enough to show
that the inequality
$$
k(d+3-k)-2\ge (2r-1)k
$$
holds for all $k\leq {\frac{d+3}{2}}$. Moreover, it is enough to
prove the latter inequality only for such natural number $k>1$
that the inequality $k(d+3-k)-2<|\bar{\Sigma}|$ holds, because
otherwise the corresponding condition on the points of the set
$\bar{\Sigma}$ is vacuous.

Now we have
$$
k(d+3-k)-2\ge k(2r-1)\iff k(r-\sum_{i=k}^{l}3(i-1)c_{i}-k)\ge
2\iff r-\sum_{i=k}^{l}3(i-1)c_{i}>k,
$$
because $k>1$. We may assume that the inequalities
$r-\sum_{i=k}^{l}3(i-1)c_{i}\le k\leq {\frac{d+3}{2}}$ and
$$
k(d+3-k)-2<|\bar{\Sigma}|
$$
hold. Let $g(x)=x(d+3-x)-2$. Then $g(x)$ is increasing for
$x<{\frac{d+3}{2}}$. Thus, we have
$$
g(k)\ge g(r-\sum_{i=k}^{l}3(i-1)c_{i}),
$$
because ${\frac{d+3}{2}}\ge k\ge r-\sum_{i=k}^{l}3(i-1)c_{i}$.
Hence, we have
$$
{\frac{(2r-1)}{3}}(r-3\sum_{i=k}^{l}ic_{i})>|\bar{\Sigma}|>g(k)\ge
(r-\sum_{i=k}^{l}3(i-1)c_{i})(2r-1)-2
$$
and
$(2r-1)(6\sum_{i=k}^{l}ic_{i}-2r)+6-9\sum_{i=k}^{l}c_{i}(2r-1)>0$.
Now $\sum_{i=k}^{l}ic_{i}<{\frac{r}{3}}$ implies
$$
(2r-1)(6\sum_{i=k}^{l}ic_{i}-2r)+6-9\sum_{i=k}^{l}c_{i}(2r-1)<6-9\sum_{i=k}^{l}c_{i}(2r-1)<6-9c_{k}(2r-1)<0,
$$
which is contradiction.
\end{proof}

Therefore, we can apply Theorem~\ref{theorem:Bese} to the blow up
of the hyperplane $\Pi$ at the points of the set
$\bar{\Sigma}\setminus\hat{P}\subset\Pi$ due to
Lemmas~\ref{lemma:double-solid-VII}, \ref{lemma:double-solid-VIII}
and \ref{lemma:double-solid-X}. The application of
Theorem~\ref{theorem:Bese} gives a curve
$C\subset\Pi\cong\mathbb{P}^{2}$ of degree
$3r-4-\sum_{i=k}^{l}3(i-1)c_{i}$ that passes trough all the points
of the set $\bar{\Sigma}\setminus\hat{P}$ and does not pass
through the point $\hat{P}=\psi(P)$. It should be pointed out that
the subset $\bar{\Sigma}\subset\Sigma^{\prime}$ may not contain
$\hat{P}\in\Sigma^{\prime}$. Namely, $\hat{P}\in\bar{\Sigma}$ if
and only if $P\in\check{\Sigma}$.

Let $G\subset\mathbb{P}^{3}$ be a cone over the curve $C$ with the
vertex $O$, where $O\in\mathbb{P}^{3}$ is the center of the
projection $\psi:\mathbb{P}^{3}\dasharrow\Pi$. Then $G$ is a
hypersurface of degree $3r-4-\sum_{i=k}^{l}3(i-1)c_{i}$ that
passes through the points of $\check{\Sigma}\setminus P$ and does
not pass through $P$. On the other hand, we already have the
hypersurface $F\subset\mathbb{P}^{3}$ of degree
$\sum_{i=k}^{l}3(i-1)c_{i}$ that passes through the points of
$\hat{\Sigma}\setminus P$ and does not pass through $P$.
Therefore, $F\cup G\subset\mathbb{P}^{3}$ is a hypersurface of
degree $3r-4$ that passes through all the points of the set
$\Sigma\setminus P$ and does not pass through the point
$P\in\Sigma$. Hence, we proved Theorem~\ref{theorem:first}.

\section{Hypersurfaces in $\mathbb{P}^{4}$.}
\label{section:hypersurfaces}

In this section we prove Theorem~\ref{theorem:second}. Let
$V\subset\mathbb{P}^{4}$ be a nodal hypersurface of degree $n$
such that $|\mathrm{Sing}(V)|\le{\frac{(n-1)^{2}}{4}}$. In order
to prove Theorem~\ref{theorem:second} it is enough to show that
the nodes of the hypersurface $V$ impose independent linear
conditions on homogeneous forms of degree $2n-5$ on
$\mathbb{P}^{4}$ due to Proposition~\ref{proposition:defect} and
Remark~\ref{remark:factoriality}. Moreover, we always may assume
that $n\ge 4$, because in the case $n\le 3$ the required claim is
trivial.

\begin{definition}
\label{definition:general-position-dva} The points of a subset
$\Gamma\subset{\mathbb P}^{r}$ satisfy the property $\bigstar$ if
at most $k(n-1)$ points of the set $\Gamma$ can lie on a curve in
$\mathbb{P}^{r}$ of degree $k\in\mathbb{N}$.
\end{definition}

Let $\Sigma=\mathrm{Sing}(V)\subset\mathbb{P}^{4}$.

\begin{proposition}
\label{proposition:nodes-in-general-position-dva} The points of
the subset $\Sigma\subset\mathbb{P}^{4}$ satisfy the property
$\bigstar$.
\end{proposition}

\begin{proof}
Let $C\subset{\mathbb P}^4$ be an irreducible and reduced curve of
degree $k$. Consider a general cone $Y\subset{\mathbb P}^4$ over
the curve $C$. Then $Y\not\subset V$ and $Y\cap V=C\cup Z$, where
$Z$ is an irreducible reduced curve of degree $k(n-1)$. Moreover,
we have the inclusion
$$
\Sigma\cap C=\mathrm{Sing}(V)\cap C\subseteq C\cap Z,
$$
but in the set-theoretic sense $|C\cap Z|\le k(n-1)$. Hence, at
most $k(n-1)$ points of the subset $\Sigma\subset\mathbb{P}^{4}$
can lie on the curve $C\subset\mathbb{P}^{4}$ of degree $k$. The
latter implies the claim.
\end{proof}

Fix a point $P\in\Sigma$. To prove that the points of
$\Sigma\subset\mathbb{P}^{4}$ impose independent linear conditions
on homogeneous forms on $\mathbb{P}^{4}$ of degree $2n-5$ it is
enough to construct a hypersurface in $\mathbb{P}^{4}$ of degree
$2n-5$ that passes through the points of the set $\Sigma\setminus
P$ and does not pass through the point $P\in\Sigma$.

\begin{lemma}
\label{lemma:hypersurface-I} Suppose that the subset
$\Sigma\subset\mathbb{P}^{4}$ is contained in some two-dimensional
linear subspace $\Pi\subset\mathbb{P}^{4}$. Then there is a
hypersurface in $\mathbb{P}^{4}$ of degree $2n-5$ that passes
through the points of the set $\Sigma\setminus P$ and does not
pass through the point $P\in\Sigma$.
\end{lemma}

\begin{proof}
Let us apply Corollary~\ref{corollary:Bese} to $\Sigma\subset\Pi$
and the number $d=2n-5\ge 3$. We must check that all the
conditions of Corollary~\ref{corollary:Bese} are satisfied. It is
clear that
$$
|\Sigma|\le {\frac{(n-1)^{2}}{4}}\Rightarrow |\Sigma|\leq {\frac{d^{2}+9d+16}{6}}%
$$
and at most $d=2n-5$ points of $\Sigma$ can lie on a line in $\Pi$
because $n\ge 4$ and the points of the subset $\Sigma\subset\Pi$
satisfy the property $\bigstar$ due to
Proposition~\ref{proposition:nodes-in-general-position-dva}.

Now we must prove that at most $k(2n-2-k)-2$ points of $\Sigma$
can lie on a curve of degree $k\leq n-1$. The case $k=1$ is
already done. Moreover, at most $k(n-1)$ points of the set
$\Sigma$ can lie on a curve of degree $k$ by
Proposition~\ref{proposition:nodes-in-general-position-dva}. Thus,
we must show that
$$
k(2n-2-k)-2\ge k(n-1)
$$
for all $k\le n-1$. Moreover, we must prove the latter inequality
only for such $k>1$ that the inequality $k(2n-2-k)-2<|\Sigma|$
holds, because otherwise the corresponding condition on the points
of the set $\Sigma$ is vacuous. Moreover,
$$
k(2n-2-k)-2\ge k(n-1)\iff n-1>k,
$$
because $k>1$. So, we may assume that $k=n-1$, but in this case
$$k(2n-2-k)-2=(n-1)^{2}>{\frac{(n-1)^{2}}{4}}\ge |\Sigma|.$$

Therefore, there is a curve $C\subset\Pi$ of degree $2n-5$ that
passes trough $\Sigma\setminus P$ and does not pass through $P$ by
Corollary~\ref{corollary:Bese}. Let $Y\subset\mathbb{P}^{4}$ be a
three-dimensional cone over $C$ with the vertex in a general line
in $\mathbb{P}^{4}$. Then $Y\subset\mathbb{P}^{4}$ is a
hypersurface of degree $2n-5$ that passes through the points of
$\Sigma\setminus P$ and does not pass through the point
$P\in\Sigma$.
\end{proof}

Fix a general two-dimensional linear subspace
$\Pi\subset\mathbb{P}^{4}$. Let $\psi:\mathbb{P}^{4}\dasharrow\Pi$
be a projection from a general line $L\subset\mathbb{P}^{4}$,
$\Sigma^{\prime}=\psi(\Sigma)$ and $\hat{P}=\psi(P)$. Then
$\psi\vert_{\Sigma}:\Sigma\to\Sigma^{\prime}$ is a bijection.

\begin{lemma}
\label{lemma:hypersurface-II} Suppose that the points in
$\Sigma^{\prime}\subset\Pi$ satisfy the property $\bigstar$. Then
there is a hypersurface in $\mathbb{P}^{4}$ of degree $2n-5$
containing $\Sigma\setminus P$ and not passing through
$P\in\Sigma$.
\end{lemma}

\begin{proof}
The proof of Lemma~\ref{lemma:hypersurface-I} implies the
existence of a curve $C\subset\Pi$ of degree $2n-5$ that passes
trough $\Sigma^{\prime}\setminus\hat{P}$ and does not pass through
$\hat{P}$. Let $Y\subset\mathbb{P}^{4}$ be a three-dimensional
cone over the curve $C$ with the vertex $L\subset\mathbb{P}^{4}$.
Then $Y\subset\mathbb{P}^{4}$ is the required hypersurface.
\end{proof}

\begin{remark}
\label{remark:projection} It seems to us that the points of the
set $\Sigma^{\prime}\subset\Pi\cong\mathbb{P}^{2}$ always satisfy
the property $\bigstar$ due to the generality in the choice of the
projection $\psi:\mathbb{P}^{4}\dasharrow\Pi$, but we fail to
prove it. In the case when
$\Sigma^{\prime}\subset\Pi\cong\mathbb{P}^{2}$ satisfy the
property $\bigstar$ the proof of Lemma~\ref{lemma:hypersurface-II}
implies a stronger result than Theorem~\ref{theorem:second}.
\end{remark}

We may assume that the points of
$\Sigma^{\prime}\subset\Pi\cong\mathbb{P}^{2}$ do not satisfy the
property $\bigstar$ and, in particular, there is a subset
$\Lambda_{k}^{1}\subset\Sigma$ such that
$|\Lambda_{k}^{1}|>k(n-1)$ and all the points of
$$
\tilde{\Lambda}_{k}^{1}=\psi(\Lambda_{k}^{1})\subset\Sigma^{\prime}\subset\Pi\cong\mathbb{P}^{2}
$$
lie on a curve $C\subset\Pi$ of degree $k$. We always may choose
$k$ to be the smallest natural number having such a property. The
latter implies that the curve $C\subset\Pi$ is irreducible and
reduced. We can iterate the construction of the subset
$\Lambda_{k}^{1}\subset\Sigma$ in the same way as in the proof of
the Theorem~\ref{theorem:first} to get disjoint subsets
$\Lambda_{j}^{i}\subset\Sigma$ for $j=k,\ldots,l\ge k$ such that
the inequality $|\Lambda_{j}^{i}|>j(n-1)$ holds, all the points of
the subset
$\tilde{\Lambda}_{j}^{i}=\psi(\Lambda_{j}^{i})\subset\Sigma^{\prime}$
lie on an irreducible reduced curve in $\Pi\cong\mathbb{P}^{2}$ of
degree $j$, and all the points in the subset
$$
\bar{\Sigma}=\Sigma^{\prime}\setminus\bigcup_{j=k}^{l}\bigcup_{i=1}^{c_{j}}\tilde{\Lambda}_{j}^{i}\subsetneq\Sigma^{\prime}\subset\Pi\cong\mathbb{P}^{2}
$$
satisfy the property $\bigstar$, where $c_{j}\ge 0$ is a number of
subsets $\tilde{\Lambda}_{j}^{i}$ and $c_{k}>0$. In particular,
$$
0\le |\bar{\Sigma}|<{\frac{(n-1)^{2}}{4}}-\sum_{i=k}^{l}c_{i}(n-1)i={\frac{n-1}{4}}(n-1-4\sum_{i=k}^{l}ic_{i}).%
$$

\begin{corollary}
\label{corollary:number-of-good-points-dva} The inequality
$\sum_{i=k}^{l}ic_{i}<{\frac{n-1}{4}}$ holds.
\end{corollary}

In particular, $\Lambda_{j}^{i}\ne\emptyset$ implies
$j<{\frac{n-1}{4}}$.

\begin{lemma}
\label{lemma:hypersurface-III} Suppose that
$\Lambda_{j}^{i}\ne\emptyset$. Let $\mathcal{M}$ be a linear
system of hypersurfaces of degree $j$ in $\mathbb{P}^{4}$ passing
through all the points in $\Lambda_{j}^{i}$. Then the base locus
of $\mathcal{M}$ is zero-dimensional.
\end{lemma}

\begin{proof}
By the construction of the set $\Lambda_{j}^{i}$ all the points of
the subset
$$
\tilde{\Lambda}_{j}^{i}=\psi(\Lambda_{j}^{i})\subset\Sigma^{\prime}\subset\Pi\cong\mathbb{P}^{2}
$$
lie on an irreducible reduced curve $C\subset\Pi$ of degree $j$.
Let $Y\subset\mathbb{P}^{4}$ be a three-dimensional cone over the
curve $C$ whose vertex is the line $L\subset\mathbb{P}^{4}$. Then
$Y$ is a hypersurface in $\mathbb{P}^{4}$ of degree $j$ that
contains all the points of the set $\Lambda_{j}^{i}$. Therefore,
$Y\in\mathcal{M}$.

Suppose that the base locus of the linear system $\mathcal{M}$
contains an irreducible reduced curve $Z\subset\mathbb{P}^{3}$.
Then $Z\subset Y$. The curves $Z$ and $C$ are irreducible and the
projection $\psi$ is sufficiently general. Therefore, $\psi(Z)=C$,
$\Lambda_{j}^{i}\subset Z$ and $\psi\vert_{Z}:Z\to C$ is a
birational morphism. In particular, $\mathrm{deg}(Z)=j$ and $Z$
contains at least $|\Lambda_{j}^{i}|>j(n-1)$ points of the subset
$\Sigma\subset\mathbb{P}^{4}$. The latter contradicts
Proposition~\ref{proposition:nodes-in-general-position-dva}.
\end{proof}

\begin{corollary}
\label{corollary:k-is-at-most-two-dva} The inequality $k\ge 2$
holds.
\end{corollary}

For every $\Lambda_{j}^{i}\ne\emptyset$ let
$\Xi_{j}^{i}\subset\mathbb{P}^{4}$ be a base locus of the linear
system of hypersurfaces of degree $j$ in $\mathbb{P}^{4}$ passing
through all the points in $\Lambda_{j}^{i}$, otherwise put
$\Xi_{j}^{i}=\emptyset$. Then $\Xi_{j}^{i}$ is a finite set by
Lemma~\ref{lemma:hypersurface-III} and
$\Lambda_{j}^{i}\subseteq\Xi_{j}^{i}$ by definition of the subset
$\Xi_{j}^{i}\subset\mathbb{P}^{4}$. Therefore, the points of the
set $\Xi_{j}^{i}\subset\mathbb{P}^4$ impose independent linear
conditions on the homogeneous forms on $\mathbb{P}^{4}$ of degree
$4(j-1)$ by Lemma~\ref{lemma:non-vanishing}. In particular, the
points of the set $\Lambda_{j}^{i}$ impose independent linear
conditions on the homogeneous forms on $\mathbb{P}^{4}$ of degree
$4(j-1)$.

\begin{lemma}
\label{lemma:hypersurface-VI} Suppose that
$\bar{\Sigma}=\emptyset$. Then there is a hypersurface in
$\mathbb{P}^{4}$ of degree $2n-5$ containing all the points of the
set $\Sigma\setminus P$ and not containing the point $P\in\Sigma$.
\end{lemma}

\begin{proof}
The points of the set $\Xi_{j}^{i}$ impose independent linear
conditions on the homogeneous forms on $\mathbb{P}^{4}$ of degree
$4(j-1)$. Therefore, for every $\Xi_{j}^{i}\ne\emptyset$
containing $P$ there is a hypersurface in $\mathbb{P}^{4}$ of
degree $4(j-1)$ that passes through the points of the set
$\Xi_{j}^{i}\setminus P$ and does not pass through the point $P$.
On the other hand, for every set $\Xi_{j}^{i}\ne\emptyset$ not
containing the point $P$ there is a hypersurface in
$\mathbb{P}^{4}$ of degree $j$ that passes through $\Xi_{j}^{i}$
and does not pass through $P$ by the definition of the set
$\Xi_{j}^{i}$. Moreover, $j<4(j-1)$, because $j\ge k\ge 2$ due to
Corollary~\ref{corollary:k-is-at-most-two-dva}. Thus, for every
non-empty set $\Xi_{j}^{i}$ there is a hypersurface
$F_{i}^{j}\subset\mathbb{P}^{4}$ of degree $3(j-1)$ that passes
through $\Xi_{j}^{i}\setminus P$ and does not pass through $P$.
Let
$$
F=\bigcup_{j=k}^{l}\bigcup_{i=1}^{c_{j}}F_{j}^{i}\subset\mathbb{P}^{4}
$$
be a possibly reducible hypersurface of degree
$\sum_{i=k}^{l}4(i-1)c_{i}$. Then $F$ passes through all the
points of the set $\Sigma\setminus P$ and does not pass through
the point $P$. Moreover, we have
$$
\mathrm{deg}(F)=\sum_{i=k}^{l}4(i-1)c_{i}<n-1\le 2n-5
$$
by Corollary~\ref{corollary:number-of-good-points-dva}, which
implies the claim.
\end{proof}

Let $\hat{\Sigma}=\cup_{j=k}^{l}\cup_{i=1}^{c_{j}}\Lambda_{j}^{i}$
and $\check{\Sigma}=\Sigma\setminus\hat{\Sigma}$. Then
$\Sigma=\hat{\Sigma}\cup\check{\Sigma}$ and
$\psi(\check{\Sigma})=\bar{\Sigma}\subset\Pi$.

\begin{remark}
\label{remark:hypersurface-through-bad-points-dva} The proof of
Lemma~\ref{lemma:hypersurface-VI} implies the existence of a
hypersurface $F\subset\mathbb{P}^{4}$ of degree
$\sum_{i=k}^{l}4(i-1)c_{i}$ that passes through all the points of
the subset $\hat{\Sigma}\setminus P\subsetneq\Sigma$ and does not
pass through the point $P\in\Sigma$.
\end{remark}

Put $d=2n-5-\sum_{i=k}^{l}4(i-1)c_{i}$. Let us check that the
subset $\bar{\Sigma}\subset\Pi\cong\mathbb{P}^{2}$ and the number
$d$ satisfy all the conditions of Theorem~\ref{theorem:Bese}. We
may assume $\hat{\Sigma}\ne\emptyset$ and
$\check{\Sigma}\ne\emptyset$.

\begin{lemma}
\label{lemma:hypersurface-VII} The inequality $d\ge 5$ holds.
\end{lemma}

\begin{proof}
We have $\sum_{i=k}^{l}4ic_{i}<n-1$ by
Corollary~\ref{corollary:number-of-good-points-dva}. Thus,
$d>n-4+4c_{k}\ge n\ge 4$.
\end{proof}

\begin{lemma}
\label{lemma:hypersurface-VIII} The inequality $|\bar{\Sigma}|\leq
{\frac{d^{2}+9d+10}{6}}$ holds.
\end{lemma}

\begin{proof}
Suppose that $|\bar{\Sigma}|>{\frac{d^{2}+9d+10}{6}}$. Then
$$
3(n-1)(n-1-4\sum_{i=k}^{l}ic_{i})>2(2n-5-\sum_{i=k}^{l}4(i-1)c_{i})^{2}+18(2n-5-\sum_{i=k}^{l}4(i-1)c_{i})+20,
$$
because
$|\bar{\Sigma}|<{\frac{n-1}{4}}(n-1-4\sum_{i=k}^{l}ic_{i})$. Let
$A=\sum_{i=k}^{l}ic_{i}$ and $B=\sum_{i=k}^{l}c_{i}$. Then
$$
3(n-1)^{2}-12(n-1)A>2(2n-1)^{2}-16A(2n-1)+32A^{2}+18(2n-1)-72A+20,
$$
because $B\ge c_{k}\ge 1$. Thus, for $n\ge 4$ we have
$$
3(n-1)^{2}>8n^{2}+28n+4+32A^{2}-A(20n+68)>
5n^{2}+12n+23>3(n-1)^{2},
$$
because $A<{\frac{n-1}{4}}$ by
Corollary~\ref{corollary:number-of-good-points-dva}.
\end{proof}

\begin{lemma}
\label{lemma:hypersurface-IX} At most $d$ of the points of the
subset $\bar{\Sigma}\subset\mathbb{P}^{2}$ lie on a line in
$\mathbb{P}^{2}$.
\end{lemma}

\begin{proof}
The points of the subset $\bar{\Sigma}\subset\mathbb{P}^{2}$
satisfy the property $\bigstar$. In particular, no more than $n-1$
of the points of $\bar{\Sigma}$ lie on a line in $\mathbb{P}^{2}$.
On the other hand, we have
$$
d=2n-5-\sum_{i=k}^{l}4(i-1)c_{i}>n-4+4c_{k}>n-1
$$
due to Corollary~\ref{corollary:number-of-good-points-dva}, which
implies the claim.
\end{proof}

\begin{lemma}
\label{lemma:hypersurface-X} At most $k(d+3-k)-2$ points of
$\bar{\Sigma}$ lie on a curve in $\mathbb{P}^{2}$ of degree $k\leq
{\frac{d+3}{2}}$.
\end{lemma}

\begin{proof}
We may assume that $k>1$ due to Lemma~\ref{lemma:hypersurface-IX}.
The points of the subset $\bar{\Sigma}\subset\mathbb{P}^{2}$
satisfy the property $\bigstar$. Thus, at most $k(n-1)$ of the
points of $\bar{\Sigma}$ lie on a curve in $\mathbb{P}^{2}$ of
degree $k$. Therefore, to conclude the proof it is enough to show
that the inequality
$$
k(d+3-k)-2\ge k(n-1)
$$
holds for all $k\leq {\frac{d+3}{2}}$. Moreover, it is enough to
prove the latter inequality only for such natural numbers $k>1$
that the inequality $k(d+3-k)-2<|\bar{\Sigma}|$ holds, because
otherwise the corresponding condition on the points of the set
$\bar{\Sigma}$ is vacuous.

Now we have
$$
k(d+3-k)-2\ge k(n-1)\iff n-1-\sum_{i=k}^{l}4(i-1)c_{i}>k,
$$
because $k>1$. Thus, we may assume that the inequalities
$$
n-1-\sum_{i=k}^{l}4(i-1)c_{i}\le k\leq {\frac{d+3}{2}}\ \ \ \text{and}\ \ k(d+3-k)-2<|\bar{\Sigma}|%
$$
hold. Let $g(x)=x(d+3-x)-2$. Then
$g(x)$ is increasing for $x<{\frac{d+3}{2}}$. Thus, we have
$$
{\frac{(n-1)}{4}}(n-1-4\sum_{i=k}^{l}ic_{i})>|\bar{\Sigma}|>g(k)\ge g(n-1-\sum_{i=k}^{l}4(i-1)c_{i}).%
$$

Let $A=\sum_{i=k}^{l}ic_{i}$ and $B=\sum_{i=k}^{l}c_{i}$. Then
the inequality
$$
{\frac{(n-1)}{4}}(n-1-4A)>4(n-1-4A+4B)(n-1)-2%
$$
holds. Therefore, we have
$$
n-1-4A>4(n-1)-16A+16B-1>4(n-1)-16A,%
$$
because $B\ge c_{k}\ge 1$. Thus, $4A>n-1$, but $A<{\frac{n-1}{4}}$
by Corollary~\ref{corollary:number-of-good-points-dva}.
\end{proof}

Therefore, we can apply Theorem~\ref{theorem:Bese} to the blow up
of the two-dimensional linear subspace $\Pi\subset\mathbb{P}^{4}$
at the points of $\bar{\Sigma}\setminus\hat{P}\subset\Pi$ due to
Lemmas~\ref{lemma:hypersurface-VII}, \ref{lemma:hypersurface-VIII}
and \ref{lemma:hypersurface-X}. The latter gives a curve
$C\subset\Pi$ of degree $2n-5-\sum_{i=k}^{l}4(i-1)c_{i}$ that
passes through all the points of the subset
$\bar{\Sigma}\setminus\hat{P}\subset\Pi\cong\mathbb{P}^{2}$ and
does not pass through the point $\hat{P}\subset\Sigma^{\prime}$.

Let $G\subset\mathbb{P}^{4}$ be a cone over the curve $C$ with the
vertex $L\subset\mathbb{P}^{4}$, where $L$ is a center of the
projection $\psi:\mathbb{P}^{4}\dasharrow\Pi$. Then
$G\subset\mathbb{P}^{4}$ is a hypersurface of degree
$2n-5-\sum_{i=k}^{l}4(i-1)c_{i}$ that passes through
$\check{\Sigma}\setminus P$ and does not pass through $P$.
However, we already have the hypersurface $F\subset\mathbb{P}^{4}$
of degree $\sum_{i=k}^{l}4(i-1)c_{i}$ that passes through
$\hat{\Sigma}\setminus P$ and does not pass through $P$.
Therefore, $F\cup G\subset\mathbb{P}^{4}$ is a hypersurface of
degree $2n-5$ that passes through $\Sigma\setminus P$ and does not
pass through $P\in\Sigma$. Thus, Theorem~\ref{theorem:second} is
proved.

\section{Calabi-Yau 3-folds.}
\label{section:Calabi-Yau}

In this section we prove Proposition~\ref{proposition:Calabi-Yau}.
Let $\pi:X\to\mathbb{P}^{3}$ be a double cover  branched over a
nodal hypersurface $S\subset\mathbb{P}^{3}$ of degree $8$ such
that $|\mathrm{Sing}(S)|\le 25$, and $V\subset\mathbb{P}^{4}$ be a
nodal hypersurface of degree $5$ such that $|\mathrm{Sing}(V)|\le
14$. Due to Proposition~\ref{proposition:defect} it is enough to
prove that the nodes of the surface $S\subset\mathbb{P}^{3}$
impose independent linear conditions on homogeneous forms of
degree $8$ on $\mathbb{P}^{3}$ and the nodes of the hypersurface
$V\subset\mathbb{P}^{4}$ impose independent linear conditions on
homogeneous forms of degree $5$ on $\mathbb{P}^{4}$.

Let $\Sigma=\mathrm{Sing}(S)\subset\mathbb{P}^{3}$ and
$\Lambda=\mathrm{Sing}(V)\subset\mathbb{P}^{4}$.

\begin{lemma}
\label{lemma:nodes-in-general-position-calabi-yau} No more than
$7k$ points of the subset $\Sigma\subset\mathbb{P}^{3}$ and no
more than $4k$ points of the subset $\Lambda\subset\mathbb{P}^{4}$
can lie on a curve of degree $k=1,2,3$.
\end{lemma}

\begin{proof}
See the proof of
Propositions~\ref{proposition:nodes-in-general-position} and
\ref{proposition:nodes-in-general-position-dva}.
\end{proof}

Fix a point $P\in\Sigma$ and a point $Q\in\Lambda$. To prove the
claim of Proposition~\ref{proposition:Calabi-Yau} we must
construct a hypersurface in $\mathbb{P}^{3}$ of degree $8$ that
passes through the points of $\Sigma\setminus P$ and does not pass
through the point $P$ and a hypersurface in $\mathbb{P}^{4}$ of
degree $5$ that passes through the points of the set
$\Lambda\setminus Q$ and does not pass through the point $Q$.

Take a general two-dimensional linear subspaces
$\Pi\subset\mathbb{P}^{3}$ and $\Omega\subset\mathbb{P}^{4}$. Let
$\psi:\mathbb{P}^{3}\dasharrow\Pi$ be a projection from a general
point $P\in\mathbb{P}^{3}$, and
$\xi:\mathbb{P}^{4}\dasharrow\Omega$ be a projection from a
general line $L\subset\mathbb{P}^{4}$. Put
$\Sigma^{\prime}=\psi(\Sigma)$, $\hat{P}=\psi(P)$,
$\Lambda^{\prime}=\xi(\Lambda)$ and $\hat{Q}=\xi(Q)$.

\begin{lemma}
\label{lemma:calabi-yau-I} No more than $7$ points of the subset
$\Sigma^{\prime}\subset\Pi$ and no more than $5$ points of the
subset $\Lambda^{\prime}\subset\Omega$ can lie on a line.
\end{lemma}

\begin{proof}
Suppose there is subset $\Theta\subset\Sigma$ such that
$|\Theta|>7$ and the points of
$\psi(\Theta)\subset\Sigma^{\prime}$ are contained in a line. Let
$\mathcal{H}$ be a linear system of hyperplanes in
$\mathbb{P}^{3}$ passing through the points of $\Theta$. Then the
base locus of $\mathcal{H}$ is zero-dimensional by
Lemma~\ref{lemma:double-solid-III}. The latter is possible only
when $|\Theta|=1$, which is a contradiction.

Suppose there is subset $\Phi\subset\Lambda$ such that
$|\Lambda|>5$ and the points of $\xi(\Phi)\subset\Lambda^{\prime}$
are contained in a line. Let $\mathcal{D}$ be a linear system of
hyperplanes in $\mathbb{P}^{4}$ passing through the points of
$\Phi$. Then the base locus of $\mathcal{D}$ is zero-dimensional
by Lemma~\ref{lemma:hypersurface-III}. The latter is possible only
when $|\Phi|=1$, which is a contradiction.
\end{proof}

\begin{lemma}
\label{lemma:calabi-yau-II} No more than $14$ points of the subset
$\Sigma^{\prime}\subset\Pi$ and no more than $10$ points of the
subset $\Lambda^{\prime}\subset\Omega$ can lie on a conic.
\end{lemma}

\begin{proof}
Suppose there is subset $\Theta\subset\Sigma$ such that
$|\Theta|>14$ and the points of $\psi(\Theta)$ are contained in a
conic $C\subset\Pi\cong\mathbb{P}^{2}$. Then $C$ is irreducible
due to Lemma~\ref{lemma:calabi-yau-I}. Let $\mathcal{H}$ be a
linear system of quadrics in $\mathbb{P}^{3}$ passing through the
points of $\Theta$. Then the base locus of the linear system
$\mathcal{H}$ is zero-dimensional by
Lemma~\ref{lemma:double-solid-III}. Take a cone
$Y\subset\mathbb{P}^{3}$ over $C$ with the vertex $O$. Then
$\Theta\subset Y$,
$\Theta\subset\mathrm{Bs}(\mathcal{H}\vert_{Y})$ and the linear
system $\mathcal{H}\vert_{Y}$ is free from base components. Let
$H_{1}$ and $H_{2}$ be general enough curves in
$\mathcal{H}\vert_{Y}$. Then $H_{i}$ is contained in the smooth
locus of the cone $Y$ and on the surface $Y$ we have
$$
8=H_{1}\cdot H_{2}\ge
\sum_{\omega\in\Theta}\mathrm{mult}_{\omega}(H_{1})\mathrm{mult}_{\omega}(H_{2})\ge
|\Theta|>14,
$$
which is a contradiction.

Let $\Phi\subset\Lambda$ be a subset such that $|\Phi|>10$.
Consider the projection $\xi$ as a composition of a projection
$\alpha:\mathbb{P}^{4}\dasharrow\mathbb{P}^{3}$ from some point
$A\in L$ and a projection $\beta:\mathbb{P}^{3}\dasharrow\Omega$
from the point $B=\alpha(L)$. The generality in the choice of the
line $L$ implies the generality of the projections $\alpha$ and
$\beta$. We claim that the points of the sets $\alpha(\Phi)$ and
$\xi(\Phi)$ do not lie on a conic in $\mathbb{P}^{3}$ and
$\Omega\cong\mathbb{P}^{2}$ respectively.

Suppose that the points of $\alpha(\Phi)$ lie on a conic
$C\subset\mathbb{P}^{3}$. Then conic $C$ is irreducible due to
Lemma~\ref{lemma:calabi-yau-I}. Let $\mathcal{D}$ be a linear
system of quadric hypersurfaces in $\mathbb{P}^{4}$ passing
through the points of $\Phi$. The proof of
Lemma~\ref{lemma:hypersurface-III} implies that the base locus of
$\mathcal{D}$ is zero-dimensional, because the points of
$\Phi\subset\mathbb{P}^{4}$ do not lie on a conic in
$\mathbb{P}^{4}$. Take a cone $W\subset\mathbb{P}^{4}$ over the
conic $C$ with the vertex $A$. Then $\Phi\subset W$. Moreover, we
have
$$
\Phi\subset\mathrm{Bs}(\mathcal{D}\vert_{W})
$$
and $\mathcal{D}\vert_{W}$ has no base components. Let $D_{1}$ and
$D_{2}$ be general curves in $\mathcal{D}\vert_{W}$. Then
$$
8=D_{1}\cdot D_{2}\ge
\sum_{\omega\in\Phi}\mathrm{mult}_{\omega}(D_{1})\mathrm{mult}_{\omega}(D_{2})\ge
|\Phi|>10,
$$
which is a contradiction. Therefore, the points of $\alpha(\Phi)$
do not lie on a conic in $\mathbb{P}^{3}$.

Suppose that the points of $\xi(\Phi)$ lie on a conic
$C\subset\Pi\cong\mathbb{P}^{3}$. Then conic $C$ is irreducible
due to Lemma~\ref{lemma:calabi-yau-I}. Let $\mathcal{B}$ be a
linear system of quadrics in $\mathbb{P}^{3}$ passing through the
points of the set $\alpha(\Phi)$. The proof of
Lemma~\ref{lemma:hypersurface-III} implies that the base locus of
$\mathcal{B}$ is zero-dimensional, because the points of
$\alpha(\Phi)$ do not lie on a conic. Take a cone
$U\subset\mathbb{P}^{3}$ over $C$ with the vertex $B$. Then
$\alpha(\Phi)\subset U$. Then
$\alpha(\Phi)\subset\mathrm{Bs}(\mathcal{B}\vert_{U})$ and the
restriction $\mathcal{B}\vert_{U}$ has no base components. Let
$B_{1}$ and $B_{2}$ be general curves in $\mathcal{B}\vert_{U}$.
Then
$$
8=B_{1}\cdot B_{2}\ge
\sum_{\omega\in\alpha(\Phi)}\mathrm{mult}_{\omega}(B_{1})\mathrm{mult}_{\omega}(B_{2})\ge
|\alpha(\Phi)|=|\Phi|>10,
$$
which is a contradiction.
\end{proof}

\begin{lemma}
\label{lemma:calabi-yau-III} There is a hypersurface in
$\mathbb{P}^{4}$ of degree $5$ that passes through the points of
the set $\Lambda\setminus Q$ and does not pass through the point
$Q\in\Lambda$.
\end{lemma}

\begin{proof}
Put $s=|\Lambda^{\prime}\setminus \hat{Q}|$. Then $s\le 13$. Let
$\pi:Y\to\Omega\cong\mathbb{P}^2$ be a blow up of the points of
the set $\Lambda^{\prime}\setminus \hat{Q}$. Then
Lemmas~\ref{lemma:calabi-yau-I} and \ref{lemma:calabi-yau-II} and
Theorem~\ref{theorem:Bese} imply the freeness of the linear system
$|\pi^{*}(\mathcal{O}_{\mathbb{P}^2}(5))-\sum_{i=1}^{s}E_{i}|$,
where $E_{i}$ is a $\pi$-exceptional curve. Let $C\subset Y$ be a
general curve in the linear system
$|\pi^{*}(\mathcal{O}_{\mathbb{P}^2}(5))-\sum_{i=1}^{s}E_{i}|$.
Then $\pi(C)\subset\Omega$ is a plane quintic curve passing
through the points of the set $\Lambda^{\prime}\setminus\hat{Q}$
and not passing through the point $Q$. The cone in
$\mathbb{P}^{4}$ over $\pi(C)$ with the vertex $L$ is the required
hypersurface.
\end{proof}

\begin{lemma}
\label{lemma:calabi-yau-IV} Suppose that at most $22$ points of
the subset $\Sigma^{\prime}\subset\Pi$ can lie on a cubic curve in
$\Pi\cong\mathbb{P}^{2}$. Then there is a hypersurface in
$\mathbb{P}^{3}$ of degree $8$ that passes through the points of
the set $\Sigma\setminus P$ and does not pass through the point
$P\in\Sigma$.
\end{lemma}

\begin{proof}
Let $\pi:Y\to\Pi\cong\mathbb{P}^{2}$ be the blow up at the points
$\{P_{1},\ldots, P_{s}\}=\Sigma^{\prime}\setminus \hat{P}$ for
$s\le 24$ and $E_{i}=\pi^{-1}(P_{i})$. Then
Lemmas~\ref{lemma:calabi-yau-I} and \ref{lemma:calabi-yau-I} and
Theorem~\ref{theorem:Bese} imply the freeness of the linear system
$|\pi^{*}(\mathcal{O}_{\mathbb{P}^2}(8))-\sum_{i=1}^{s}E_{i}|$.
Let $C\in
|\pi^{*}(\mathcal{O}_{\mathbb{P}^2}(8))-\sum_{i=1}^{k}E_{i}|$ be a
general enough curve. Then $\pi(C)\subset\Pi$ is a plane octic
curve that passes through the points of
$\Sigma^{\prime}\setminus\hat{P}$ and does not pass through the
point $\hat{P}$. Hence, the cone in $\mathbb{P}^{3}$ over the
curve $\pi(C)\subset\Pi$ with the vertex $O$ is the required
hypersurface.
\end{proof}

\begin{lemma}
\label{lemma:calabi-yau-V} Suppose that there is a subset
$\Upsilon\subset\Sigma$ such that $|\Upsilon|>22$ and all the
points of the set $\psi(\Upsilon)$ lie on a cubic curve in
$\Pi\cong\mathbb{P}^{2}$. Then there is a hypersurface in
$\mathbb{P}^{3}$ of degree $8$ that passes through the points of
$\Sigma\setminus P$ and does not pass through the point $P$.
\end{lemma}

\begin{proof}
Let $\mathcal{H}$ be a linear system of cubic hypersurfaces in
$\mathbb{P}^{3}$ passing through the points of the set $\Upsilon$.
Then the base locus of $\mathcal{H}$ is zero-dimensional by
Lemma~\ref{lemma:double-solid-III}.

Suppose $P\in\Upsilon$. Then there is a hypersurface
$F\subset\mathbb{P}^{3}$ of degree $6$  that passes through the
points of $\Upsilon\setminus P$ and does not pass through the
point $P$ by Lemma~\ref{lemma:non-vanishing}. On the other hand,
the subset $\Sigma\setminus\Upsilon\subset\mathbb{P}^{3}$ contains
at most $2$ points. Hence, there is a quadric
$G\subset\mathbb{P}^{3}$ that passes through the points of
$\Sigma\setminus\Upsilon$ and does not pass through $P$. Thus,
$F\cup G$ is the required hypersurface.

In the case when $P\not\in\Upsilon$ and
$P\in\mathrm{Bs}(\mathcal{H})$ we can repeat every step of the
proof of the previous case. In the case when $P\not\in\Upsilon$
and $P\not\in\mathrm{Bs}(\mathcal{H})$ there is a cubic
hypersurface in $\mathbb{P}^{3}$  that passes through the points
of $\Upsilon$ and does not pass through the point $P$, which
easily implies the existence of the required hypersurface.
\end{proof}

Hence, Proposition~\ref{proposition:Calabi-Yau} is proved. It
seems to us that the bounds for nodes in
Proposition~\ref{proposition:Calabi-Yau} can be improved using
the methods of \cite{Rei88}, \cite{EinLa93} and \cite{Ka97}
instead of Theorem~\ref{theorem:Bese}.

\section{Non-isolated singularities.}
\label{section:applications}

In this section we prove Theorems~\ref{theorem:third} and
\ref{theorem:forth}. Let $\tau:U\to\mathbb{P}^{s}$ be a double
cover branched over a hypersurface $F$ of degree $2r$ such that
$D_{1}\cap\cdots\cap D_{s-3}$ is a $\mathbb{Q}$-factorial nodal
3-fold, where $D_{i}$ is a general divisor in
$|\tau^{*}(\mathcal{O}_{\mathbb{P}^{s}}(1))|$. Let
$W\subset\mathbb{P}^{r}$ be a hypersurface of degree $n$ such that
$H_{1}\cap\cdots\cap H_{r-4}$ is a $\mathbb{Q}$-factorial nodal
3-fold, where $H_{i}$ is a general hyperplane section of $W$. We
may assume that $s\ge 4$ and $r\ge 5$. We must show that the group
$\mathrm{Cl}(U)$ is generated by $D_{1}$ and the group
$\mathrm{Cl}(W)$ is generated by $H_{1}$.

\begin{remark}
\label{remark:normality} The varieties $U$ and $W$ are normal (see
Proposition~8.23 in \cite{Ha77}).
\end{remark}

Let $D$ be a general divisor in
$|\tau^{*}(\mathcal{O}_{\mathbb{P}^{s}}(1))|$ and $H$ be a general
hyperplane section of $W$.

\begin{lemma}
\label{lemma:simple-vanishing-I} The groups
$H^{1}(\mathcal{O}_{U}(-nD))$ and $H^{1}(\mathcal{O}_{W}(-nH))$
vanish for every $n>0$.
\end{lemma}

\begin{proof}
In the case when the varieties $U$ and $W$ have mild singularities
the claim is implies by the Kawamata-Viehweg vanishing (see
\cite{Ka82}, \cite{Vi82}). In general let us prove the claim by
induction on $s$ and $r$. We consider only the vanishing of
$H^{1}(\mathcal{O}_{U}(-nD))$, because the proof of the vanishing
of the cohomology group  $H^{1}(\mathcal{O}_{W}(-nH))$ is
identical.

Suppose that $s=4$. Then we have an exact sequence of sheaves
$$
0\to\mathcal{O}_{U}(-(n+1)D)\to\mathcal{O}_{U}(-nD)\to\mathcal{O}_{D}(-nD)\to 0%
$$
for any $n\in\mathbb{Z}$. Therefore, we have an exact sequence of
the cohomology groups
$$
0\to H^{1}(\mathcal{O}_{U}(-(n+1)D))\to H^{1}(\mathcal{O}_{U}(-nD))\to H^{1}(\mathcal{O}_{D}(-nD))\to\cdots%
$$
for $n>0$. However, the 3-fold $D$ is nodal by assumption. Thus,
the group $H^{1}(\mathcal{O}_{D}(-nD))$ vanishes by the
Kawamata-Viehweg vanishing. Hence, we have
$$
H^{1}(\mathcal{O}_{U}(-D))\cong H^{1}(\mathcal{O}_{U}(-2D))\cong\cdots\cong H^{1}(\mathcal{O}_{U}(-nD))%
$$
for every $n>0$. On the other hand, the group
$H^{1}(\mathcal{O}_{U}(-nD))$ vanishes for $n\gg 0$ by the lemma
of Enriques-Severi-Zariski (see \cite{Za52} or Corollary~7.8 in
\cite{Ha77}).

Suppose that $s>4$. Then we have an exact sequence of sheaves
$$
0\to\mathcal{O}_{U}(-(n+1)D)\to\mathcal{O}_{U}(-nD)\to\mathcal{O}_{D}(-nD)\to 0%
$$
for any $n\in\mathbb{N}$. Therefore, we have an exact sequence of
the cohomology groups
$$
0\to H^{1}(\mathcal{O}_{U}(-(n+1)D))\to H^{1}(\mathcal{O}_{U}(-nD))\to H^{1}(\mathcal{O}_{D}(-nD))\to\cdots%
$$
for $n>0$. However, the group $H^{1}(\mathcal{O}_{D}(-nD))$
vanishes by the induction. Hence,
$$
H^{1}(\mathcal{O}_{U}(-D))\cong
H^{1}(\mathcal{O}_{U}(-2D))\cong\cdots\cong
H^{1}(\mathcal{O}_{U}(-nD))%
$$
for $n>0$, but $H^{1}(\mathcal{O}_{U}(-nD))=0$ for $n\gg 0$ by the
lemma of Enriques-Severi-Zariski.
\end{proof}

Consider a Weil divisor $G$ on $U$. Let us prove by the induction
on $s$ that $G\sim kD$ for some $k\in\mathbb{Z}$. Suppose that
$s=4$. Then the 3-fold $D$ is nodal and $\mathbb{Q}$-factorial by
assumption. Moreover, the group $\mathrm{Cl}(D)$ is generated by
the class of the divisor $R\vert_{D}$ due to
Remark~\ref{remark:factoriality}, where $R$ is a general divisor
in the linear system $|D|$. Thus, there is an integer $k$ such
that $G\vert_{D}\sim kR\vert_{D}$. Let $\Delta=G-kR$. Then the
sequence of sheaves
$$
0\to\mathcal{O}_{U}(\Delta)\otimes \mathcal{O}_{U}(-D)\to\mathcal{O}_{U}(\Delta)\to\mathcal{O}_{D}\to 0%
$$
is exact, because $\mathcal{O}_{U}(\Delta)$ is locally free in the
neighborhood of $D$. Hence, the sequence
$$
0\to H^{0}(\mathcal{O}_{U}(\Delta))\to H^{0}(\mathcal{O}_{D})\to H^{1}(\mathcal{O}_{U}(\Delta)\otimes \mathcal{O}_{U}(-D))%
$$
is exact.

\begin{lemma}
\label{lemma:simple-vanishing-II} The group
$H^{1}(\mathcal{O}_{U}(\Delta)\otimes \mathcal{O}_{U}(-nD))$
vanishes for every $n>0$.
\end{lemma}

\begin{proof}
The sheaf $\mathcal{O}_{U}(\Delta)$ is reflexive (see
\cite{Har80}). Thus, there is an exact sequence of sheaves
$$
0\to\mathcal{O}_{U}(\Delta)\to\mathcal{E}\to\mathcal{F}\to 0%
$$
where $\mathcal{E}$ is a locally free sheaf and $\mathcal{F}$ is a
torsion free sheaf. Hence, the sequence of groups
$$
H^{0}(\mathcal{F}\otimes\mathcal{O}_{U}(-nD))\to H^{1}(\mathcal{O}_{D}(\mathrm{\Delta})\otimes\mathcal{O}_{D}(-nD))\to H^{1}(\mathcal{E}\otimes\mathcal{O}_{U}(-nD))%
$$
is exact. However, for $n\gg 0$ the cohomology group
$H^{0}(\mathcal{F}\otimes\mathcal{O}_{U}(-nD))$ vanishes because
the sheaf $\mathcal{F}$ is torsion free, and the cohomology group
$H^{1}(\mathcal{E}\otimes\mathcal{O}_{U}(-nD))$ vanishes by the
lemma of Enriques-Severi-Zariski (see \cite{Za52} or Corollary~7.8
in \cite{Ha77}). Therefore, the cohomology group
$H^{1}(\mathcal{O}_{U}(\Delta)\otimes\mathcal{O}_{U}(-nD))$
vanishes for $n\gg 0$.

Now consider an exact sequence of sheaves
$$
0\to\mathcal{O}_{U}(\Delta)\otimes\mathcal{O}_{U}(-(n+1)D)\to\mathcal{O}_{U}(\Delta)\otimes\mathcal{O}_{U}(-nD)\to\mathcal{O}_{D}(-nD)\to 0%
$$
and the induced sequence of the cohomology groups
$$
0\to H^{1}(\mathcal{O}_{U}(\Delta)\otimes\mathcal{O}_{U}(-(n+1)D))\to H^{1}(\mathcal{O}_{U}(\Delta)\otimes\mathcal{O}_{U}(-nD))\to H^{1}(\mathcal{O}_{D}(-nD))\to\cdots%
$$
for $n>0$. Then the group $H^{1}(\mathcal{O}_{D}(-nD))$ vanishes
by Lemma~\ref{lemma:simple-vanishing-I}. Hence, we have
$$
H^{1}(\mathcal{O}_{U}(\Delta)\otimes\mathcal{O}_{U}(-D))\cong H^{1}(\mathcal{O}_{U}(\Delta)\otimes\mathcal{O}_{U}(-2D))\cong\cdots\cong H^{1}(\mathcal{O}_{U}(\Delta)\otimes\mathcal{O}_{U}(-nD))%
$$
for $n>0$, but we already proved that
$H^{1}(\mathcal{O}_{U}(-nD))$ vanishes for $n\gg 0$.
\end{proof}

Therefore, $H^{0}(\mathcal{O}_{U}(\Delta))\cong\mathbb{C}$.
Similarly $H^{0}(\mathcal{O}_{U}(-\Delta))\cong\mathbb{C}$. Thus,
the Weil divisor $\Delta$ is rationally equivalent to zero and
$G\sim kD$ in the case $s=4$.

Suppose that $s>4$. By the induction we may assume that the group
$\mathrm{Cl}(D)$ is generated by the class of the divisor
$R\vert_{D}$, where $R$ is a general divisor in $|D|$. Thus, there
is an integer $k$ such that $G\vert_{D}\sim kR\vert_{D}$. Put
$\Delta=G-kR$. Then the sequence of sheaves
$$
0\to\mathcal{O}_{U}(\Delta)\otimes \mathcal{O}_{U}(-D))\to\mathcal{O}_{U}(\Delta)\to\mathcal{O}_{D}\to 0%
$$
is exact, because $\mathcal{O}_{U}(\Delta)$ is locally free in the
neighborhood of $D$. Therefore, the sequence
$$
0\to H^{0}(\mathcal{O}_{U}(\Delta))\to H^{0}(\mathcal{O}_{D})\to H^{1}(\mathcal{O}_{U}(\Delta)\otimes \mathcal{O}_{U}(-D))%
$$
is exact. However, the proof of the
Lemma~\ref{lemma:simple-vanishing-II} holds for $s>4$. Thus, the
cohomology group $H^{1}(\mathcal{O}_{U}(\Delta)\otimes
\mathcal{O}_{U}(-D))$ vanishes. Hence,
$H^{0}(\mathcal{O}_{U}(\Delta))\cong\mathbb{C}$. Same arguments
prove that $H^{0}(\mathcal{O}_{U}(-\Delta))\cong\mathbb{C}$.
Therefore, the Weil divisor $\Delta$ is rationally equivalent to
zero and $G\sim kD$. Thus, we proved Theorem~\ref{theorem:third}.
We omit the proof of Theorem~\ref{theorem:forth}, because it is
identical to the proof of the Theorem~\ref{theorem:third}.

\section{Birational rigidity.}
\label{section:birational-rigidity}

In this section we prove
Proposition~\ref{proposition:fourfold-double-cover}. Let $\xi:Y\to
\mathbb{P}^{4}$ be a double cover branched over a hypersurface
$F\subset\mathbb{P}^{4}$ of degree $8$ such that the hypersurface
$F$ is smooth outside of a smooth curve $C\subset F$, the
singularity of the hypersurface $F$ in a sufficiently general
point of the curve $C$ is locally isomorphic to the singularity
$$
x_{1}^{2}+x_{2}^{2}+x_{3}^{2}=0\subset\mathbb{C}^{4}\cong\mathrm{Spec}(\mathbb{C}[x_{1},x_{2},x_{3},x_{4}]),
$$
the singularities of $F$ in other points of $C$ are locally
isomorphic to the singularity
$$
x_{1}^{2}+x_{2}^{2}+x_{3}^{2}x_{4}=0\subset\mathbb{C}^{4}\cong\mathrm{Spec}(\mathbb{C}[x_{1},x_{2},x_{3},x_{4}]),
$$
and a general 3-fold in $|-K_{Y}|$ is $\mathbb{Q}$-factorial. Then
$Y$ is a Fano 4-fold with terminal singularities and
$-K_{Y}\sim\xi^{*}(\mathcal{O}_{\mathbb{P}^{4}}(1))$. Moreover,
$\mathrm{Cl}(Y)$ and $\mathrm{Pic}(Y)$ are generated by the
divisor $-K_{Y}$ by Theorem~\ref{theorem:third}. Hence, $Y$ is a
Mori fibration (see \cite{KMM}). We must prove that the 4-fold $Y$
is a unique Mori fibration birational to $Y$ and
$\mathrm{Bir}(Y)=\mathrm{Aut}(Y)$. It is well known that the
latter implies the finiteness of the group $\mathrm{Bir}(Y)$.

Suppose that either $Y$ is not birationally rigid or
$\mathrm{Bir}(Y)\ne\mathrm{Aut}(Y)$. Then
Theorem~\ref{theorem:Nother-Fano} imply the existence of a linear
system $\mathcal{M}$ on $Y$ such that $\mathcal{M}$ has no fixed
components and the singularities of $(X,
{\frac{1}{n}}\mathcal{M})$ are not canonical, where
$\mathcal{M}\sim -nK_{Y}$. Thus, there is a rational number
$\mu<{\frac{1}{n}}$ such that $(X, \mu\mathcal{M})$ is not
canonical, i.e. $\mathbb{CS}(Y, \mu\mathcal{M})\ne\emptyset$.

Let $Z$ be an element of the set $\mathbb{CS}(Y,\mu\mathcal{M})$.
Then $\mathrm{mult}_{Z}(\mathcal{M})>n$.

\begin{lemma}
\label{lemma:smooth-points} The subvariety $Z\subset Y$ is not a
smooth point of $Y$.
\end{lemma}

\begin{proof}
Suppose that $Z$ is a smooth point of $Y$. Then
$$
\mathrm{mult}_{Z}(\mathcal{M}^{2})>4n^{2}
$$
by Theorem~\ref{theorem:Iskovskikh}. Take general divisors $H_{1}$
and $H_{2}$ in $|-K_{Y}|$ containing $Z$. Then
$$
2n^{2}=\mathcal{M}^{2}\cdot H_{1}\cdot H_{2}\ge \mathrm{mult}_{Z}(\mathcal{M}^{2})\mathrm{mult}_{Z}(H_{1})\mathrm{mult}_{Z}(H_{2})>4n^{2}%
$$
which is a contradiction.
\end{proof}

\begin{lemma}
\label{lemma:singular-points} The subvariety $Z\subset Y$ is not a
singular point of $Y$.
\end{lemma}

\begin{proof}
Let $\xi(Z)=O$. Then $O$ is a singular point of the hypersurface
$F\subset\mathbb{P}^{4}$. Therefore, the point $O$ is contained in
the curve $C\subset F$ by assumption. There are two possible
cases, i.e. either the singularity of $F$ in the point $O$ is
locally isomorphic to the singularity
$$x_{1}^{2}+x_{2}^{2}+x_{3}^{2}=0\subset\mathbb{C}^{4}\cong\mathrm{Spec}(\mathbb{C}[x_{1},x_{2},x_{3},x_{4}]),$$
or the singularity of $F$ in the point $O$ is locally isomorphic
to the singularity
$$x_{1}^{2}+x_{2}^{2}+x_{3}^{2}x_{4}=0\subset\mathbb{C}^{4}\cong\mathrm{Spec}(\mathbb{C}[x_{1},x_{2},x_{3},x_{4}]),$$
where $x_{1}=x_{2}=x_{3}$ are local equations of the curve
$C\subset F$. Let us call the former case ordinary and the latter
case non-ordinary.

Let $X$ be a sufficiently general divisor in the linear system
$|-K_{Y}|$ passing through the point $Z$. Then the double cover
$\xi$ induces the double cover $\tau:X\to\mathbb{P}^{3}$ ramified
in an octic surface. The singularities of $X\setminus Z$ are
ordinary double points. Moreover, $Z$ is an ordinary double point
of $X$ in the ordinary case. In the non-ordinary case the
singularity of the 3-fold $X$ at the point $Z$ is locally
isomorphic to
$$x_{1}^{2}+x_{2}^{2}+x_{3}^{2}+x_{4}^{3}=0\subset\mathbb{C}^{4}\cong\mathrm{Spec}(\mathbb{C}[x_{1},x_{2},x_{3},x_{4}]).$$

Let $\mathcal{D}=\mathcal{M}\vert_{X}$ and $H=-K_{Y}\vert_{X}$.
Then the linear system $\mathcal{D}$ has no fixed components and
$\mathcal{D}\sim nH$. Moreover, $Z\in\mathbb{LCS}(X,
\mu\mathcal{D})$ by Theorem~\ref{theorem:log-adjunction}. In
particular, $Z\in\mathbb{CS}(X, \mu\mathcal{D})$.

Let $f:V\to X$ be a blow up of $Z$, $E=f^{-1}(Z)$ and
$\mathcal{H}$ be a proper transform of the linear system
$\mathcal{D}$ on $V$. Then $V$ is smooth in the neighborhood of
$E$ and $E$ is isomorphic to a quadric surface in
$\mathbb{P}^{3}$. In the ordinary case $E$ is smooth. In the
non-ordinary case the quadric surface $E$ has one singular point
$P\in E$, i.e. the surface $E$ is isomorphic to a quadric cone in
$\mathbb{P}^{3}$. Note, that $K_{V}\sim E$.

Let $\mathrm{mult}_{Z}(\mathcal{D})\in\mathbb{N}$ such that
$\mathcal{H}\sim f^{*}(nH)-\mathrm{mult}_{Z}(\mathcal{D})E$. Then
$\mathrm{mult}_{Z}(\mathcal{D})>n$ in the ordinary case by
Theorem~\ref{theorem:Corti}. On the other hand, in the
non-ordinary case we have the inequality
$\mathrm{mult}_{Z}(\mathcal{D})>{\frac{n}{2}}$ due to
Proposition~\ref{proposition:non-simple-double-point}.

By construction the linear system $|f^{*}(H)-E|$ is free and gives
a morphism $\psi:V\to\mathbb{P}^{2}$ such that
$\psi=\phi\circ\tau\circ f$, where
$\phi:\mathbb{P}^{3}\dasharrow\mathbb{P}^{2}$ is a projection from
the point $O$. Moreover, the restriction
$\psi\vert_{E}:E\to\mathbb{P}^{2}$ is a double cover. Let $L$ be a
sufficiently general fiber of the morphism $\psi$. Then $L$ is a
smooth curve of genus $2$ and $L\cdot E=L\cdot f^{*}(H)=2$. Thus,
$$
L\cdot \mathcal{H}=L\cdot
f^{*}(nH)-\mathrm{mult}_{Z}(\mathcal{D})L\cdot
E=2n-2\mathrm{mult}_{Z}(\mathcal{D})\ge 0,
$$
because $\mathcal{H}$ has no base components. Hence,
$\mathrm{mult}_{Z}(\mathcal{D})\le n$. In particular, the ordinary
case is impossible and it remains to eliminate the non-ordinary
case.

The inequalities $\mathrm{mult}_{Z}(\mathcal{D})\le n$ and
$\mu<{\frac{1}{n}}$, the equivalence
$$
K_{V}+\mu\mathcal{H}\sim
f^{*}(K_{X}+\mu\mathcal{D})+(1-\mu\mathrm{mult}_{Z}(\mathcal{D}))E
$$
and $Z\in\mathbb{CS}(X, \mu\mathcal{D})$ imply the existence of a
proper irreducible subvariety $S\subset E$ such that
$S\in\mathbb{CS}(V,
\mu\mathcal{H}+(\mu\mathrm{mult}_{Z}(\mathcal{D})-1)E)$. In
particular, $S\in\mathbb{CS}(V, \mu\mathcal{H})$.

Suppose that $S$ is a curve. Then
$\mathrm{mult}_{S}(\mathcal{H})>n$. Let $L_{\omega}$ be a fiber of
$\psi$ passing through a general point $\omega\in S$. Then
$L_{\omega}$ spans a divisor in $V$ when we vary $\omega$ on $C$.
Hence,
$$
L\cdot \mathcal{H}=L\cdot
f^{*}(nH)-\mathrm{mult}_{Z}(\mathcal{D})L\cdot
E=2n-2\mathrm{mult}_{Z}(\mathcal{D})\ge
\mathrm{mult}_{\omega}(L_{\omega})\mathrm{mult}_{S}(\mathcal{H})>n,
$$
which contradicts the inequality
$\mathrm{mult}_{Z}(\mathcal{D})>{\frac{n}{2}}$.

Therefore, $S$ is a point on $E$. Then
$\mathrm{mult}_{S}(\mathcal{H})>n$ and
$\mathrm{mult}_{S}(\mathcal{H}^{2})>4n^{2}$ by
Theorem~\ref{theorem:Iskovskikh}, because $S$ is smooth on $V$. It
is easy to see that the point $S$ is not a vertex $P$ of the
quadric cone $E$, because the numerical intersection of a general
ruling of $E$ with a general divisor in $\mathcal{H}$ is equal to
$\mathrm{mult}_{Z}(\mathcal{D})\le n$. Let $\Gamma$ be a fiber of
the morphism $\psi$ that passes through the point $S$ and $D$ be a
general divisor in the linear system $|f^{*}(H)-E|$ that passes
through the point $S$. Then $\Gamma\subset D$. Note, that $\Gamma$
may be reducible and singular, but we always have the inequality
$\mathrm{mult}_{S}(\Gamma)\le 2$, because $\tau\circ f(\Gamma)$ is
a line passing through the point $O$ and $\tau\vert_{f(\Gamma)}$
is a double cover.

Suppose that $\Gamma$ is irreducible. Let
$$
\mathcal{H}^{2}=\lambda\Gamma+T,
$$
where $\lambda\in\mathbb{Q}$ and $T$ is a one-cycle such that
$\Gamma\not\subset\mathrm{Supp}(T)$. Then the inequalities
$$
\mathrm{mult}_{S}(T)>4n^{2}-\lambda\mathrm{mult}_{S}(\Gamma)\ge 4n^{2}-2\lambda%
$$
hold. On the other hand, the inequalities
$$
\mathrm{mult}_{S}(T)\le \mathrm{mult}_{S}(T)\mathrm{mult}_{S}(D)\le T\cdot D=\mathcal{H}^{2}\cdot D=2n^{2}-\mathrm{mult}^{2}_{Z}(\mathcal{D})<{\frac{7}{4}}n^{2}%
$$
holds. Thus, we have $\lambda>{\frac{9}{8}}n^{2}$. Let $\tilde{D}$
be a general divisor in $|f^{*}(H)|$. Then
$$
2n^{2}=\tilde{D}\cdot\mathcal{H}^{2}\ge\lambda\Gamma\cdot\tilde{D}=2\lambda>{\frac{9}{4}}n^{2},%
$$
which is a contradiction.

Therefore, the fiber $\Gamma$ is reducible. Then
$\Gamma=\Gamma_{1}\cup\Gamma_{2}$, where $\Gamma_{i}$ is a smooth
rational curve such that $\tau\circ f(\Gamma_{1})=\tau\circ
f(\Gamma_{2})$ is a line in $\mathbb{P}^{3}$ containing point $O$.
Let
$$
\mathcal{H}^{2}=\lambda_{1}\Gamma_{1}+\lambda_{2}\Gamma_{2}+T,
$$
where $\lambda_{i}\in\mathbb{Q}$ and $T$ is a one-cycle such that
$\Gamma_{i}\not\subset\mathrm{Supp}(T)$. Then the inequalities
$$
{\frac{7}{4}}n^{2}>2n^{2}-\mathrm{mult}^{2}_{Z}(\mathcal{D})\ge T\cdot D\ge\mathrm{mult}_{S}(T)>4n^{2}-\lambda_{1}-\lambda_{2}%
$$
hold. Thus, $\lambda_{1}+\lambda_{2}>{\frac{9}{4}}n^{2}$. Hence,
we have
$$
2n^{2}=\tilde{D}\cdot\mathcal{H}^{2}\ge\lambda_{1}\Gamma_{1}\cdot\tilde{D}+\lambda_{2}\Gamma_{2}\cdot\tilde{D}=\lambda_{1}+\lambda_{2}>{\frac{9}{4}}n^{2}%
$$
for a general divisor $\tilde{D}\in |f^{*}(H)|$, which is a
contradiction.
\end{proof}

\begin{lemma}
\label{lemma:codimension-big} The subvariety $Z\subset Y$ is not a
curve.
\end{lemma}

\begin{proof}
Suppose $Z$ is a curve. Let $X$ be a  general divisor in
$|-K_{Y}|$ and $P$ be a point in the intersection $Z\cap X$. Then
$X$ is a nodal Calabi-Yau 3-fold. The point $P$ is smooth on the
3-fold $X$ if and only if $Z\not\subset\mathrm{Sing}(X)$. In the
case $Z\subset\mathrm{Sing}(X)$ the point $P$ is an ordinary
double point on $X$. Moreover, $P\in\mathbb{CS}(X,
\mu\mathcal{D})$, where $\mathcal{D}=\mathcal{M}\vert_{X}$. In the
case when the point $P$ is smooth on $X$ we can proceed as in the
proof of Lemma~\ref{lemma:smooth-points} to get a contradiction.
In the case when the point $P$ is an ordinary double point on $X$
we can proceed as in the proof of
Lemma~\ref{lemma:singular-points} to get a contradiction.
\end{proof}

\begin{lemma}
\label{lemma:codimension-two} The subvariety $Z\subset Y$ is not a
surface.
\end{lemma}

\begin{proof}
Suppose $Z$ is a surface. Then $\mathrm{mult}_{Z}(\mathcal{M})>n$.
Let $V$ be a general divisor in the linear system $|-K_{Y}|$,
$S=Z\cap V$ and $\mathcal{D}=\mathcal{M}\vert_{V}$. Then $V$ is a
nodal Calabi-Yau 3-fold, the linear system $\mathcal{D}$ has no
base components, $S\subset V$ is a an irreducible reduced curve
and $\mathrm{mult}_{S}(\mathcal{D})>n$. The double cover $\xi$
induces a double cover $\tau:V\to\mathbb{P}^{3}$ ramified in a
nodal hypersurface $G\subset\mathbb{P}^{3}$ of degree $8$.

Take a sufficiently general divisor $H$ in
$|\tau^{*}(\mathcal{O}_{\mathbb{P}^{3}}(1))|$. Then
$$
2n^{2}=\mathcal{D}^{2}\cdot
H\ge\mathrm{mult}^{2}_{S}(\mathcal{D})S\cdot H> n^{2}S\cdot H,%
$$
which implies $S\cdot H=1$. Hence, $\tau(S)$ is a line in
$\mathbb{P}^{3}$ and $\tau\vert_{S}$ is an isomorphism.

Suppose that $\tau(S)\not\subset G$. Then there is a smooth
rational curve $\tilde{S}\subset V$ such that $S\ne\tilde{S}$ and
$\tau(S)=\tau(\tilde{S})$. Take a sufficiently general surface
$D\in |\tau^{*}(\mathcal{O}_{\mathbb{P}^{3}}(1))|$ passing through
the curve $S$. Then $D$ is smooth outside of $S\cap\tilde{S}$.
Moreover, the surface $D$ is smooth in every point of
$S\cap\tilde{S}$ that is smooth on $V$, and $D$ has an ordinary
double point in every point of $S\cap\tilde{S}$ that is an
ordinary double point on $V$. On the other hand, at most $4$ nodes
of the hypersurface $G\subset\mathbb{P}^{3}$ can lie on the line
$\tau(S)$, i.e. $|\mathrm{Sing}(D)|\le 4$. The sub-adjunction
formula (see \cite{KMM}, \cite{Ko91}) implies
$$
(K_{D}+\tilde{S})\vert_{\tilde{S}}=K_{\tilde{S}}+\mathrm{Diff}_{\tilde{S}}(0)%
$$
and $\mathrm{deg}(\mathrm{Diff}_{\tilde{S}}(0))={\frac{k}{2}}$,
where $k=|\mathrm{Sing}(D)|$. Thus, the self-intersection
$\tilde{S}^{2}$ is negative on the surface $D$, because
$K_{D}\cdot\tilde{S}=1$. Put $\mathcal{H}=\mathcal{D}\vert_{D}$. A
priori the linear system $\mathcal{H}$ can have a base component.
However, the generality in the choice of $D$ implies
$$
\mathcal{H}=\mathrm{mult}_{S}(\mathcal{D})S+\mathrm{mult}_{\tilde{S}}(\mathcal{D})\tilde{S}+\mathcal{B}%
$$
where $\mathcal{B}$ is a linear system on $D$ having no base
components. Moreover, the equivalence
$$
(n-\mathrm{mult}_{\tilde{S}}(\mathcal{D}))\tilde{S}\sim_{\mathbb{Q}}(\mathrm{mult}_{S}(\mathcal{D})-n)S+\mathcal{B}%
$$
and $\tilde{S}^{2}<0$ imply
$\mathrm{mult}_{\tilde{S}}(\mathcal{D})>n$. Take a general divisor
$H$ in $|\tau^{*}(\mathcal{O}_{\mathbb{P}^{3}}(1))|$. Then
$$
2n^{2}=\mathcal{D}^{2}\cdot
H\ge\mathrm{mult}^{2}_{S}(\mathcal{D})S\cdot H+\mathrm{mult}^{2}_{\tilde{S}}(\mathcal{D})\tilde{S}\cdot H>n^{2}S\cdot H+n^{2}\tilde{S}\cdot H=2n^{2},%
$$
which is a contradiction.

Therefore, we have $\tau(S)\subset G$. Let $O$ be a general point
on $\tau(S)$ and $\Pi$ be a hyperplane in $\mathbb{P}^{3}$ that
tangents $G$ at the point $O$. Consider a sufficiently general
line $L\subset\Pi$ passing through $O$. Let $\hat{L}=\tau^{-1}(L)$
and $\hat{O}=\tau^{-1}(O)$. Then $\hat{L}$ is singular at
$\hat{O}$. Therefore, the curve $\hat{L}$ is contained in the base
locus of the linear system $\mathcal{D}$, because otherwise
$$
2n=\hat{L}\cdot\mathcal{D}\ge\mathrm{mult}_{\hat{O}}(\hat{L})\mathrm{mult}_{\hat{O}}(\mathcal{D})\ge 2\mathrm{mult}_{S}(\mathcal{D})>2n%
$$
which is impossible. On the other hand, the curve $\hat{L}$ spans
a divisor in $V$ when we vary the line $L$ in $\Pi$. The latter is
impossible, because $\mathcal{D}$ has no base components.
\end{proof}

Therefore, Proposition~\ref{proposition:fourfold-double-cover} is
proved.

\end{document}